%% file: POMS_Template_V1.tex
\documentclass[poms,final,nonblindrev]{poms1_V1} 

\OneAndAHalfSpacedXI 
\RequirePackage{etex}
\RequirePackage{fix-cm}

 \def\ntpss{\textsc{NTP}\big([a_j]_{j \in \mathcal{A}}, W\big)}
\def\feasS{\mathscr{F}_{{NTP} }}
\def\optS{\mathscr{O}_{{NTP} }}
\def\probname{\textsc{NTP}}
\def\detprobname{\textsc{NTP}}
\def\dettwoprobname{\textsc{NTP2}}
\def\stochprobname{\textsc{DNTP}}

\usepackage{graphicx}
\usepackage{subfigure, epsfig}
\usepackage{natbib}
\usepackage{amsmath}
\usepackage{lmodern}
\usepackage{anyfontsize}
\usepackage{subcaption}
\usepackage{pgfgantt}
\usepackage[mathscr]{eucal}
\usepackage{adjustbox}
\usepackage{multicol}
\usepackage{amssymb}
\usepackage{multirow}
\usepackage{varwidth}
\usepackage{hyperref}
\usepackage{xcolor}
\usepackage{hanging}
\usepackage{linegoal}
\usepackage{bbm}
\usepackage[skip=5pt,font=scriptsize]{caption}
\usepackage{lipsum}
\usepackage[ruled,linesnumbered]{algorithm2e}
\usepackage[nameinlink]{cleveref}

 \bibpunct[, ]{(}{)}{,}{a}{}{,}%
 \def\bibfont{\small}%
\newcommand{\quotes}[1]{``#1''}
\TheoremsNumberedThrough     
\ECRepeatTheorems

\EquationsNumberedThrough    

\begin{document}


\RUNAUTHOR{Gawas, Legrain and Rousseau}

\RUNTITLE{NTP for On-Demand Personnel Scheduling}

\TITLE{Notification Timing for On-Demand Personnel Scheduling}

\ARTICLEAUTHORS{%
\AUTHOR{Prakash Gawas }
\AFF{ Polytechnique Montreal,  GERAD,  CIRRELT, and HANALOG, Montréal, Canada,  \EMAIL{prakash.gawas@polymtl.ca}} 
\AUTHOR{  Antoine Legrain}
\AFF{Polytechnique Montreal,  GERAD,  CIRRELT, and HANALOG, Montréal, Canada, \EMAIL{antoine.legrain@polymtl.ca}} 
\AUTHOR{  Louis-Martin Rousseau}
\AFF{Polytechnique Montreal,  CIRRELT, and HANALOG,  Montréal, Canada, \EMAIL{louis-martin.rousseau@polymtl.ca}} 

} 

\ABSTRACT{%
Modern business models have introduced service systems that can tap into a large pool of casual employees with flexible working hours who are paid based on piece rates to fill in for on-demand work. Such systems have been successfully implemented in various domains like ride-sharing, deliveries, and microtasks. However, since these casual employees have infrequent engagement and may not have sufficient experience, the quality of work is a significant issue, especially in service industries.  We present a novel casual personnel scheduling system designed to provide experienced casual employees to service companies to optimize their operations with a dynamic, data-driven approach. Like traditional on-call systems, it contacts casual personnel in order of seniority to inform them about the availability of on-demand work.  However, flexibility in the system allows employees time to think and freely choose shifts available to them. It also allows senior employees to replace (bump) junior employees from the current schedule if no other shift of their choice is available, provided certain conditions are met. The replacements, though allowed by the system, cause inconvenience and dissatisfaction among employees. The management seeks to quickly assign all available shifts to employees while ensuring minimum bumps. However, there is uncertainty associated with the time at which an employee would select a shift.  The decision faced by the management is to determine when to notify an employee to avoid bumps.  We first show that this problem is $\mathcal{NP}-complete$ even in the case of perfect information. We then propose a two-stage stochastic formulation for the dynamic problem and develop a heuristic algorithm that approximates optimal policy using a threshold structure. These policies are calibrated using offline solutions where all uncertainty is pre-known, enabling us to fine-tune the policy structure. Testing on real-world data shows that our policy outperforms the current strategy used by our industry partner.
}%


\maketitle

%


\input{intro}

\input{pdes}
\input{dntp}

\input{exp}
\input{conclusion}



\bibliographystyle{pomsref}

 \let\oldbibliography\thebibliography
 \renewcommand{\thebibliography}[1]{%
    \oldbibliography{#1}%
    \baselineskip14pt 
    \setlength{\itemsep}{10pt}
 }



%
%
%

\ECSwitch 

\ECHead{Notification Timing for On-Demand Personnel Scheduling - E-Companion}


\input{appendix}

\end{document}

%% file: intro.tex
\section{Introduction}\label{introduction}

The Personnel (or shift) scheduling problem was first introduced in the 1950s through the pioneering works of \cite{dantzig1954comment} and \cite{edie1954traffic}. Since then, it has been the subject of extensive research (e.g. \cite{de2015workforce}, \cite{van2013personnel}). Traditionally, service industries plan their regular shift schedules weeks in advance, often relying on perceived cost-saving measures such as overtime or part-time employment to address uncertainties that emerge later. However, these conventional scheduling approaches have become less effective in the current environment, characterized by unpredictable demand, seasonal fluctuations, growing complexity, and shifting employee work preferences. The lack of a modern, adaptive scheduling system can lead to significant challenges for both management and employees, including reduced service quality, operational inefficiencies, and lower profitability.
\subsection{On-Demand Labor Platforms}

The gig economy (\cite{graham2019gig}) has recently emerged as a new business model that has experienced significant growth over the past decade (\cite{donovan2016does}), offering an alternative approach for industries. It operates as an on-demand labor platform, where workers take on short-term tasks or freelance projects rather than traditional full-time employment. These workers, commonly referred to as freelancers, casual workers, contingent workers, gig workers, or independent contractors, engage with platforms such as ride-sharing services like Uber and Lyft, food delivery services like GrubHub, Seamless, and InstaCart, customer contact centers like Liveops, home services like Handy, and online platforms for web development and data analysis like Upwork and Fiverr. These for-profit companies function as on-demand platforms (\cite{allon2023impact}, \cite{musalem2023balancing}), helping businesses meet short-term labor needs by connecting them to a flexible, crowd-sourced workforce. A defining feature of these platforms is their ability to tap into a large, crowd-sourced pool of workers to manage demand surges while controlling labor costs through variable pay instead of fixed wages. 

On-demand labor platforms operate through four key elements: (1) the task, (2) the client posting the task, (3) the worker with the relevant skills, and (4) the platform that manages the interaction. Workers, often freelancers, create profiles highlighting their skills and rates, while clients post task listings specifying the project scope, budget, and deadlines. These tasks can be open to all workers or limited to specific individuals, allowing workers to then choose the ones that suit their preferences. The flexible nature of these platforms lets workers set their schedules and switch between platforms, which is a hallmark of the gig economy. It offers companies the advantage of hiring on-demand labor and paying only for completed tasks.

Despite these advantages, online labor platforms face several domain-specific challenges in task design, pricing models, and task allocation (\cite{slivkins2014online}, \cite{bhatti2020general}). Addressing all these concerns in a unified model is complex, leading different platforms to adopt varied approaches. Consider task allocation (\cite{karachiwalla2021understanding}, \cite{zhen2021crowdsourcing}), which will be the primary focus of this paper. On sites like Fiverr and Upwork, workers actively search for projects with varying scopes and requirements, submitting applications and bids—either proposing a fixed project price or an hourly rate. Even when clients set compensation, workers still retain some control over which tasks they pursue. For instance, on Amazon Mechanical Turk (MTurk), workers browse micro-tasks, each with a predetermined piecework rate set by the client. In contrast, location-based platforms often rely on algorithms to assign tasks, with pay rates set by the platform itself. For example, Uber or Lyft drivers receive ride requests with only seconds to decide, often without knowing the passenger's destination. Similarly, TaskRabbit workers, known as "taskers", have a limited time window—ranging from minutes to hours—to respond to job offers, depending on urgency and client preferences.

\subsection{On-Demand Shift Scheduling System}

In this paper, we study an on-demand system service companies use to perform various daily tasks. Examples include amusement parks, airports, security services, and call centers. These organizations must dynamically adjust their workforce based on real-time demand to balance labor costs from underutilization with the need to respond to customer demands, thereby generating revenue. Tasks are organized into shifts that require specific skills to be completed within set timeframes. For instance, in an amusement park (client), shifts might involve roles like ticket checking, cleaning, or waiting staff—positions that typically require low skill levels and offer modest pay. The client posts shifts on the platform, which acts as a self-scheduling system, informing workers of available on-demand opportunities. Workers can accept, decline, or ignore these offers, giving them a high degree of flexibility. Our industrial partner, Merinio, currently operates such a platform, managing approximately 15,000 employees for 1.4 million shifts annually, showcasing the system's scalability and efficiency. Similar platforms include Wonolo, GigSmart, Instawork, and Shiftgig.

Clients using on-demand platforms aim to deliver high-quality customer service by maintaining optimal staffing levels with an efficient and flexible workforce. However, managing this with casual workers presents challenges, such as limited screening, inconsistent engagement, and lower wages. Additionally, some casual workers may struggle with motivation and communication. On most platforms, tasks are open to the entire workforce, resulting in a high number of applicants. However, clients typically favor workers with the highest reputations, awarding them the most work opportunities (\cite{yoganarasimhan2013value}), as they are viewed as more reliable and capable of providing quality service. In contrast, workers with lower skills, poor reputations, or those new to the platform often struggle to secure assignments. Reputation is generally based on client reviews and points rating.

The self-scheduling on-demand system we examine offers flexibility by gradually opening job opportunities based on employee seniority. Initially, jobs are made available to the most senior employees, and if unclaimed, they become accessible to junior employees. Workers are notified as soon as a job is made available to them. Senior workers, viewed as more reliable, are given the advantage of claiming jobs earlier and, under certain conditions, can even claim shifts already chosen by junior employees. While this system helps maintain service quality, it can create challenges for junior employees, who may struggle to secure their preferred shifts or, at times, any work. This introduces the challenge of optimizing job opening notifications, which we call the Notification Timing Problem. This paper aims to develop an optimal notification policy that balances the inconvenience to junior employees with company staffing goals. Unlike traditional on-call systems, which lack flexibility and can result in irregular schedules and health issues (\citet{golden2015irregular}, \citet{nicol2004call}), the on-demand system integrates employee preferences, improving job satisfaction and reducing turnover. Real-time platforms also ease the administrative burden for employers to manage and select workers from a large workforce.

\subsection{Research Contribution} \label{rc}
In particular, we make the following contributions:
\begin{enumerate}
    \item We present a novel and significant problem related to timing job opening notifications to employees for last-minute shift assignments in Personnel Scheduling, specifically within dynamic on-demand labor platforms. To our knowledge, this is the first study to focus on optimizing the operations of a dynamic on-demand scheduling system, where employees can decide whether they want to work in response to the notifications sent. The objective of the problem is to strike a tradeoff between quick shift assignment and employee inconvenience caused by the on-demand system.
    \item We then show that even a simple and offline version of the problem with complete information is $\mathcal{NP}-complete$ in \Cref{ntp}. We show a reduction from the subset sum problem by giving examples, intuition, and proof.
    \item We consider the (real) stochastic and dynamic version of the same problem and show that any dynamic policy can be as bad as possible compared to its offline counterpart. We develop an operating policy based on deterministic offline solutions obtained with simulated data. The methodology aggregates the offline solutions to get a nonparametric policy function.
    \item   We compare the performance with heuristic policy functions and demonstrate that offline solution-based policy outperforms them on real-world data from our industrial partner. We extended our results with different employee shift preference distributions and showed that the developed offline policy remains better.
\end{enumerate}

%% file: pdes.tex
\section{Problem Description and Computational Complexity}\label{ntp}

In this section, we formally define the Notification Timing Problem and highlight its significance. We then consider a simplified offline variant, introducing its notations and presenting a Mixed Integer Programming (MIP) formulation. The main result in this section shows that the offline version is $\mathcal{NP}$-complete through a reduction, supported by intuitive reasoning and illustrative examples.

\subsection{Problem Description}\label{syswork}

The client company initially assesses the number of on-demand shifts required for a specific upcoming day based on a precise forecast of future demand. This required number of shifts is then fed into an electronic system overseen by the platform. The system generates a list of eligible employees according to their skills and availability, and an operations manager initiates the communication process on the client's behalf. The electronic system notifies workers about available shifts, monitors their responses, and updates the shift schedule as employees accept the assignments.

Notifications of job openings are sent to employees on their phones in the order of seniority, which is determined by experience and reputation. The number of notifications sent at any time follows a preset operating policy. Employees can take their time before responding to available shifts, as they might be busy, need to consider their options or wish to consult with family. If employees do not answer, it indicates their unavailability for that day. The interval between the system's notification and the employee's response is called \quotes{response delay}. Once an employee accepts a shift, that shift is marked as occupied. As more employees respond and claim shifts, overall shift occupancy increases. Employees generally have different preferences that influence their choices for shifts.

The system also allows senior employees to select shifts that junior employees have already claimed, provided they respond within a specified time frame known as the \quotes{cutoff}. This results in a schedule change where the originally scheduled employee is \quotes{bumped}. Senior employees have this privilege because they are the most experienced and productive. The cutoff period begins for an employee from the moment a notification is sent ($s_i$) to that employee. The employee $i$ then responds back at time $e_i = s_i + r_i$. If employees respond after the cutoff, they can only choose from unoccupied shifts, thereby losing their seniority advantage. The system encourages prompt shift selection by rewarding employees who consistently respond quickly with a higher seniority level, giving them access to a broader range of shifts. Bumped employees are immediately notified and can choose from the remaining available shifts at that time. This also means that they can bump some junior employees who hold their next preferred shift, and thus, one response can initiate multiple bumps. The scheduling of shifts is finalized either when all employees have responded or when the planning period ends.

The uncertainty in the system arises from the fact that the employee response delay is not known in advance. If too many employees are contacted early in the scheduling horizon, this can lead to multiple bumps, especially if some senior employees take longer to make decisions. Frequent disruptions or bumps can cause frustration among junior employees and may even lead them to leave the system, which is particularly problematic given the transient nature of the workforce. Employee retention is, therefore, crucial. Conversely, if too few employees are contacted to avoid bumps, there may not be enough time for responses within the planning horizon, potentially resulting in unassigned shifts.

This creates a challenge for the system in determining how many notifications are to be sent at any given time, considering the current state of the system. From a management perspective, it is essential to establish a notification policy that minimizes bumps while ensuring all shifts are filled. However, these two goals often conflict, adding complexity to the problem. For example, to reduce bumps, an ideal policy might be not to send any notifications, but this would leave shifts unfilled. Conversely, if the priority is to ensure that all shifts are covered, notifying everyone at the outset might seem optimal, but it could lead to excessive bumps.

We refer to this problem of determining the optimal timing for contacting employees as the Notification Timing Problem ({\probname}). The decision involves finding a notification policy that balances these conflicting objectives. Specifically, at any given time, the policy should assess the system's current state and determine when to contact the next employee. Therefore, the {\probname} is an online sequential decision problem with uncertainty, where each decision incurs a cost in terms of either bumps or unassigned shifts.

\subsection{Offline  \textsc{NTP} Formulation}\label{prelim}
We define a simplified offline version of the (\textsc{NTP}) with specific assumptions and simplifications. Although this version differs slightly from the real-time, online problem, it highlights the complexity of solving it, even under simplified conditions. In this offline(full information) problem, all employee response delays and shift preferences are known at the outset of the planning process. Thus one can theoretically solve the problem and have policy ready beforehand.  All employee response delays are considered to be shorter than the planning horizon. While, in reality, some response times may exceed this horizon, such employees are excluded from consideration in the offline version since they would be unable to select shifts. Response delays are rounded to the nearest integer. \Cref{tab:all_nots} outlines the notations used in this problem.

\begin{table}[ht!]
\caption{Model, Parameters, and Variables of the \textsc{NTP}}
\label{tab:all_nots}
\parbox{1\linewidth}{
\centering
\resizebox{\linewidth}{!}{
  {\def\arraystretch{1.25}  
\begin{tabular}{llll}
\hline
\textbf{Sets}& &&\\
\hline
$\mathcal{E}  $ & &  &set of employees in the order of seniority with one being most senior \\
$\mathcal{L}$ & & & set of shifts \\
$\mathcal{X}_i = \{l^1_i, l^2_i,\dots, l^M_i\} $ & & & ordered set of employee preferences, $l^1_i \in \mathcal{L}$ \\
\hline
\textbf{Parameters}&&&\\
\hline
$M$ & & &number of employees\\
$L$ & & &number of shifts\\
$\mathbf{r} = [r_i]_{i \in \mathcal{E}}$ & & & response delays for all employees \\
$H$ && &planning horizon\\
\hline
\textbf{Variables}&&&\\
\hline
$s_i$ & $i \in \mathcal{E}$ & $s_i \geq 0 $ & time when a job opening notification is sent to  employee $i$ \\
$e_i$ & $i \in \mathcal{E}$ &  $e_i \geq 0$ &time of response for employee $i$ \\
$y_{ij}$ & $i,j \in \mathcal{E}$ & $y_{ij} \in \{0, 1\} $ & 1 if $i$ bumps $j$ else 0 \\
$b_i = \sum\limits_{j \in \mathcal{E}: i <j}y_{ij}$ & $i \in \mathcal{E}$ & $b_i \geq 0 $ &number bumps caused by employee $i$ \\
\hline

\end{tabular}}}

}
\end{table}

The planning horizon is discretized into units, and notifications are sent only at these discrete time steps. The planning horizon $H$ is much greater than the number of employees available after discretization. We assume that there is no cutoff time—senior employees can respond and bump junior employees at any point during the planning process. We consider that there are the same number of shifts ($L = M$) to fill as there are employees, and all employees are eligible to work all shifts.

In actual online instances, the number of employees typically exceeds the number of available shifts. Some employees may not want to work and may not respond. On the other hand, the total number of employees willing to work can exceed the number of shifts. Since this information is known in advance, one can easily filter employees in such cases to ensure $L = M$. The number of shifts that remain vacant in online instances will depend on how many employees can respond within the horizon. Given that we know the response delays in advance for offline instances, we can efficiently compute the number of vacant shifts beforehand. Thus, this cost becomes constant, and the objective is simply to minimize total bumps, given the knowledge of vacancies in the offline formulation.

Each employee is identified by their seniority number $i \in \mathcal{E} = \{1, 2, \dots, M\}$. For any two employees $i$, $j$ in $\mathcal{E}$, if $i < j$, employee $i$ is more senior than employee $j$, and thus must be notified before or simultaneously with employee $j$. Let $L_i \in \mathcal{X}_i$ represent the next available shift for employee $i$ at the current time $t$. Initially, at $t = 0$, $L_i = l^1_i$, where $l^1_i$ is employee $i$'s most preferred shift from their preference set, $\mathcal{X}_i = \{l^1_i, l^2_i, \dots, l^M_i\}$. Additionally, let $I_l$ denote the employee currently assigned to shift $l \in \mathcal{L}$ at time $t$. At $t = 0$, we have $I_l = \phi$, indicating that all shifts are unoccupied at the start.

When employee $i$ responds, the following steps are followed:
\begin{itemize}
    \item[1.] \textbf{Selection of Shift}: Employee $i$ will always select their most preferred available shift, $L_i$. If this shift is currently unoccupied $(I_{L_i} = \phi)$, it is directly assigned to employee $i$, meaning $I_{L_i} = i$. This shift is then no longer available for any employee $j$ such that $i < j$. Employee $i$'s next preferred available shift $L_i$ is updated, and is set to the next shift in their preference order, as long as that shift is occupied by a junior employee or still unoccupied.
    \item[2.] \textbf{Bumping Scenario}: If $L_i$ is currently occupied by a junior employee $j_1$ (i.e., $I_{L_i} = j_1$ and $i < j_1$), a bump occurs. In this case, employee $i$ replaces $j_1$ in shift $L_i$, and we update $I_{L_i} = i$. $j_1$ now must be reassigned to another shift. The process is repeated for employee $j_1$, where steps 1 and 2 are applied to find a new shift for them. If $j_1$'s next preferred shift is also occupied, further bumps may occur.
    \item [3.] \textbf{Availability of Shifts}: At any point in the simulation, if there is no shift available for some junior employee $j$ as a result of bumping, then $L^i = \phi$, referring to a null shift. This will happen only when senior employees to $j$ have selected all available shifts.  
\end{itemize}   

 Hence, by design of the system, a response by an employee $i$ may result in a chain of bumps. Let $\mathcal{B}_i$ refer to a chain of employees in the order of seniority who are bumped due to $i$ responding at time $t$. We define the total bumps by employee $i$ at time $t$ as the total schedule changes due to a response, which is $|\mathcal{B}_i| = b_i$. We can construct the bump chain $\mathcal{B}_i$ as follows, given $i$ has responded at time $t$. A pseudo-code Algorithm 3 to compute the chain is also given in the E-Companion.
 
\begin{align}
\mathcal{B}_i = \{j_1,\dots, j_k: I_{L_i} = j_1, I_{L_{j_k}} = \phi, I_{L_{j_o}}= j_{o+1}, {j_o} < j_{o+1}, 1 \leq o < k\} 
\label{bumpset}
\end{align}

In this case, employee $j_1$ is directly bumped by employee $i$, while other employees are indirectly bumped in a chain as $j_1$ is reassigned to a new shift, followed by $j_2$, and so on. The final employee in the chain, $j_k$, either takes an unoccupied shift or is assigned to a null shift, indicating that they are no longer able to work. Due to discretization, it is possible that multiple employees may respond at the same time. In such instances, we compute $\mathcal{B}_i$ for each responding employee in order of seniority, ensuring no bumps occur between employees responding within the same time frame. This can result in a situation where a junior employee might be counted as being bumped twice, even though there is only one actual shift change for that employee at that time. This reflects the real-time nature of operations, where the system is continuously updated. While two senior employees responding at the exact same moment may seem possible, in practice, real-time updates prevent this scenario. The system logs the response of one senior employee before the other, causing the junior employee to experience two bumps sequentially, even though both senior employees responded in quick succession.

\subsection{Employee Preferences}
Employee preferences are difficult to model, and our industrial partner does not store any data on employees. This data is difficult to collect as asking employees to rank their preferences every time is not feasible. Even if this data is available, all employee preferences are expected to change daily. Further, learning a model over such a dataset would induce errors. Hence, in the absence of the preference data, our aim in solving the problem is to obtain a schedule to send notifications to minimize all \textbf{potential} instances of bumps. A potential instance of a bump is when a senior employee responds after a junior employee. It is equivalent to saying we want to minimize the total count where a senior employee responds after a junior employee. We then show that {potential} instances of bumps are a tight upper bound on the total bumps given any set of employee preferences. Given an employee $i$ has responded at time $e_i$, the number of potential bumps $p_i$ by employee $i$ equals $|\mathcal{P}_i|$  where 

\begin{align}
\mathcal{P}_i = \{j: j \in \mathcal{E}, e_i > e_j, i < j \}  
\label{pbumpset}
\end{align}

Let $S = [(s_i, e_i)]_{i \in \mathcal{E}}$ denote a notification schedule. $S$ is a feasible schedule for the instance if all the employees are notified in the order of seniority ($s_i \leq s_{i+1}$) and all of them respond within the horizon ($e_i \leq H$).   Also, let $B_S = \sum_{i \in \mathcal{E}} b_i$ denote the total bumps suffered by the employees when following the notification schedule $S$, where each employee $i$ induces $b_i = |\mathcal{B}_i| = \sum_{j \in \mathcal{E}} y_{ij}$ bumps. The variable $b_i$ will be used later in the paper. Nevertheless, we define it here. Then, the following theorem holds.

\begin{theorem}
\label{theo3}
Given an instance $I = \{H, M, \mathbf{r}\}$, let $S^*_\mathcal{X}$  denote an optimal schedule with preferences $\mathcal{X}$ such that $\mathcal{X} = \{{\mathcal{X}_i}\}_{i \in \mathcal{E}}$. Also, let $S^*_p$ be the schedule that minimizes the total potential bumps for this instance, and  $S^*_I$ be the schedule that minimizes the total bumps with identical preferences. Then, $B_{S^*_\mathcal{X}} \leq B_{S^*_p} = B_{S^*_I} $.
\end{theorem}

\begin{corollary}
\label{cor:0}
Minimizing the total potential bumps is equivalent to reducing total bumps under the same employee preferences.
\end{corollary}

The proof of \Cref{theo3} is provided in the E-Companion. Consequently, the {\detprobname} with the same preferences and equal employees and shifts represents the worst-case scenario regarding bumps. Diverse employee preferences serve only to decrease the occurrence of bumps. Henceforth, {\detprobname} will specifically denote the problem in which the total potential bumps are minimized. All variables related to bumps, $b_i, y_{ij}, B_S$ will refer to potential bumps. Any reference to a bump from this point on will correspond to a potential bump unless mentioned otherwise. In the following sections, we show that this problem is hard; hence, minimizing bumps under some preferences is also hard.

\subsection{MIP Model}
The MIP formulation for the simplified  {\detprobname} is given by  (\eqref{eq1} - \eqref{eq5}).  This formulation only decides when to notify employees. It does not explicitly allow shifts to employees. To have a feasible schedule, we simply need to ensure that all employees respond within the time horizon $H$. Note that in this formulation, we merely minimize the total number of bumps in the final schedule. The vacancy cost is not included as it is easy to know precisely how many shifts one can schedule using the offline information indicated earlier.

\resizebox{16cm}{!}{
\centering
\makebox[\textwidth]{\begin{minipage}{\dimexpr\textwidth-4\fboxsep-8\fboxrule\relax}
\begin{equation}
{\textbf{\textsc{MIP}}_{\textbf{NTP}}} := \quad \min \quad \sum_{i\in \mathcal{E}}\sum_{j\in \mathcal{E}:i<j} y_{ij} \label{eq1}
\end{equation}
\begin{align}
\textbf{s.t.} \quad\quad\quad\quad &s_i   \leq s_j & \forall i,j\in \mathcal{E}, i < j \label{eq2}\\
&e_i \leq  H & \forall i\in \mathcal{E} \label{eq3}\\
&e_i  = s_i + r_i & \forall i\in \mathcal{E} \label{eq4}\\
& e_i - e_j \leq (r_i - r_j)y_{ij} & \forall i,j\in \mathcal{E}, i < j, r_i > r_j \label{eq5}\\
&y_{ij} \in \{0,1\} & \forall i,j\in \mathcal{E}, i < j \label{eq6}\\
&s_{i}, e_i \geq 0 & \forall i\in \mathcal{E} \label{eq7}
\end{align}
\end{minipage}}
}
\vspace{0.5cm}

The objective function minimizes total bumps. \eqref{eq2} ensures the feasibility of the schedule so that employees are notified according to seniority. Constraint \eqref{eq3} makes sure that all employees respond before the horizon. The start and end times should match the response delay through \eqref{eq4}.   Finally, Constraint \eqref{eq5} counts each potential bump instance in the schedule. A potential bump occurs when any employee $i$ responds after any junior employee $j$ that is $e_i > e_j$.  This equation represents the $Big-M$ constraint. $(r_i - r_j)$ is a valid value of $Big-M$ because first $i < j$ and the seniority constraint must be respected between $i$ and $j$, that is, $s_i \leq s_j$. Then, the bound of $(r_i - r_j)$ is only met when $s_i = s_j$, when $i, j$ are called at once.  In addition, we only include the constraint when $i < j$ and $ r_i > r_j$. When $ r_i \leq r_j$, $i$ cannot bump $j$ due to the seniority constraint as $i$ would respond earlier than $j$. This constraint also counts bumps as they would be in real-time when there are multiple responses simultaneously.  Note that in this formulation, we do not explicitly track potential bumps. To count the potential total bumps an employee suffers, one just needs to count how many senior employees have responded after that employee.  This also implies that employees who get bumped are allocated the next shift in the identical preference order.

\subsection{Preliminaries}

Let {\detprobname($I$)}  denote  the {\detprobname}  defined by instance $I = \{H, M, \mathbf{r}\}$.   Also, let $\feasS(I)$ denote the set of all feasible solutions for {\detprobname($I$)}. Given any schedule $S $, the \emph{makespan}, denoted by $C(S)$, is defined as the time required to schedule all shifts. 
\begin{equation*}
    C(S) = \max_{i \in \mathcal{E}}\{e_i\}
\end{equation*}
The makespan is feasible if $C(S) \leq H$. Consider any two employees $ i, j \in \mathcal{E}$ where $i <j$. A No Bump Schedule (NBS) is any schedule $S$ such that $B(S) = 0$, i.e., there cannot be two employees $i, j$, and $i< j$ such that $e_i > e_{j}$. Note that this schedule will not have potential bumps or realized bumps. Any NBS, if feasible for a given instance, is also optimal. Let $C^*_0 $ denote the minimum makespan of an NBS such that $B(S) = 0$. 

The following proposition gives an expression to calculate $C^*_0 $ for the  {\detprobname}.
\begin{proposition}
\label{prop_makespan_r}
    The minimum makespan of a NBS is $C^*_0 = r_1 + \sum\limits_{i=1}^{M-1} {(r_{i+1} - r_{i})^+}$.
\end{proposition} 
The following example illustrates the problem.

\paragraph{\textbf{Example}:}  Consider an instance $I$ with $M = 6$ employees, a planning horizon $H = 10$, and response times $\mathbf{r} = [4, 1, 5, 3, 2, 5]$. \Cref{gantt1}(a) presents an NBS for this instance, with a makespan of 11, indicated by the blue line. For clarity, we refer to each employee as $E_i$. An NBS is feasible for {\detprobname($I$)} if the condition $H \geq C^*_0$ is satisfied. However, since $H = 10$, it is evident that an NBS is not feasible in this case. Let $\optS(I) \subset \feasS(I)$ represent the set of optimal schedules that minimize potential bumps for instance $I$. Multiple optimal solutions can exist for the same instance: \Cref{gantt1}(b) and (c) show two such schedules, $S_1, S_2\in \optS(I)$, respectively. Introducing potential bumps into the schedule can reduce the total time needed to complete all shifts within the planning horizon. This is demonstrated in \Cref{gantt1}(a) and (b). In \Cref{gantt1}(a), employee $E_2$ is notified at time 3, while in \Cref{gantt1}(b), $E_2$ is notified one unit earlier, at time 2. Notifying $E_2$ earlier introduces a bump for that employee, but it allows the system to notify all subsequent employees one unit earlier, speeding up the process. \Cref{gantt1}(c) illustrates a similar situation for employee $E_5$.

This example highlights why introducing bumps is sometimes necessary to complete a schedule within the given horizon.

\begin{figure}[!ht]
\caption{(a) NBS : $\mathbf{C_0^* = 11}$,\hfill (b) $\mathbf{S_1^*}$ : $\mathbf{B(S_1^*) = 1}$ and $\mathbf{C(S_1^*) = 10}$,\hfill (c) $\mathbf{S_2^*}$ : $\mathbf{B(S_2^*) = 1}$ and $\mathbf{C(S_2^*) = 10}$ }
\label{gantt1}
\begin{minipage}{0.33\textwidth}
\begin{ganttchart}[
x unit  = 0.35cm,
hgrid = true,
vgrid=dashed, bar/.append style={fill=gray!50},
y unit chart=11,
vrule/.style={very thick, blue},
title height=.4,
title height=0.8,
y unit title=0.6cm,
inline,
link/.style={-latex, draw=red, fill=red}
]{0}{11}

\gantttitle{\textbf{Notification Schedule}}{12} \\
\gantttitlelist[title label font=\tiny,]{1,...,12}{1}\\
\ganttgroup[inline=false,group label font=\scriptsize, group/.style={draw=none, fill=none}]   {   E1}{  0}{ 11   }
\ganttbar[bar label font=\tiny\color{black},]{4}{0}{3} \\
\ganttgroup[inline=false,group label font=\scriptsize, group/.style={draw=none, fill=none}]   {   E2}{  0}{ 11   }
\ganttbar[bar label font=\tiny\color{black},]{1}{3}{3}\\
\ganttgroup[inline=false,group label font=\scriptsize, group/.style={draw=none, fill=none}]   {   E3}{  0}{ 11   }
\ganttbar[bar label font=\tiny\color{black},]{5}{3}{7}\\
\ganttgroup[inline=false,group label font=\scriptsize, group/.style={draw=none, fill=none}]   {   E4}{  0}{ 11   }
\ganttbar[bar label font=\tiny\color{black},]{3}{5}{7}\\
\ganttgroup[inline=false,group label font=\scriptsize, group/.style={draw=none, fill=none}]   {   E5}{  0}{ 11   }
\ganttbar[bar label font=\tiny\color{black},]{2}{6}{7}\\
\ganttgroup[inline=false,group label font=\scriptsize, group/.style={draw=none, fill=none}]   {   E6}{  0}{ 11   }
\ganttbar[bar label font=\tiny\color{black},]{5}{6}{10}\\
\ganttvrule[
vrule/.append style={red}]{$H = 10$}{9}
\ganttvrule[
vrule/.append style={blue}]{}{10}
\end{ganttchart}
\end{minipage}\hfill
\begin{minipage}{0.33\textwidth}
\begin{ganttchart}[
x unit  = 0.35cm,
hgrid = true,
vgrid=dashed, bar/.append style={fill=gray!50},
y unit chart=11,
vrule/.style={very thick, blue},
title height=.4,
title height=0.8,
y unit title=0.6cm,
inline,
link/.style={-latex, draw=red, fill=red}
]{0}{11}

\gantttitle{\textbf{Notification Schedule}}{12} \\
\gantttitlelist[title label font=\tiny,]{1,...,12}{1}\\
\ganttgroup[inline=false,group label font=\scriptsize, group/.style={draw=none, fill=none}]   {   E1}{  0}{ 11   }
\ganttbar[bar label font=\tiny\color{black},]{4}{0}{3} \\
\ganttgroup[inline=false,group label font=\scriptsize, group/.style={draw=none, fill=none}]   {   E2}{  0}{ 11   }
\ganttbar[bar label font=\tiny\color{black},]{1}{2}{2}\\
\ganttgroup[inline=false,group label font=\scriptsize, group/.style={draw=none, fill=none}]   {   E3}{  0}{ 11   }
\ganttbar[bar label font=\tiny\color{black},]{5}{2}{6}\\
\ganttgroup[inline=false,group label font=\scriptsize, group/.style={draw=none, fill=none}]   {   E4}{  0}{ 11   }
\ganttbar[bar label font=\tiny\color{black},]{3}{4}{6}\\
\ganttgroup[inline=false,group label font=\scriptsize, group/.style={draw=none, fill=none}]   {   E5}{  0}{ 11   }
\ganttbar[bar label font=\tiny\color{black},]{2}{5}{6}\\
\ganttgroup[inline=false,group label font=\scriptsize, group/.style={draw=none, fill=none}]   {   E6}{  0}{ 11   }
\ganttbar[bar label font=\tiny\color{black},]{5}{5}{9}\\
\ganttvrule[
vrule/.append style={blue}]{$H = 10$}{9}
\end{ganttchart}

\end{minipage}\hfill
\begin{minipage}{0.33\textwidth}
\begin{ganttchart}[
x unit  = 0.35cm,
hgrid = true,
vgrid=dashed, bar/.append style={fill=gray!50},
y unit chart=11,
vrule/.style={very thick, blue},
title height=.4,
title height=0.8,
y unit title=0.6cm,
inline,
link/.style={-latex, draw=red, fill=red}
]{0}{11}

\gantttitle{\textbf{Notification Schedule}}{12} \\
\gantttitlelist[title label font=\tiny,]{1,...,12}{1}\\
\ganttgroup[inline=false,group label font=\scriptsize, group/.style={draw=none, fill=none}]   {   E1}{  0}{ 11   }
\ganttbar[bar label font=\tiny\color{black},]{4}{0}{3} \\
\ganttgroup[inline=false,group label font=\scriptsize, group/.style={draw=none, fill=none}]   {   E2}{  0}{ 11   }
\ganttbar[bar label font=\tiny\color{black},]{1}{3}{3}\\
\ganttgroup[inline=false,group label font=\scriptsize, group/.style={draw=none, fill=none}]   {   E3}{  0}{ 11   }
\ganttbar[bar label font=\tiny\color{black},]{5}{3}{7}\\
\ganttgroup[inline=false,group label font=\scriptsize, group/.style={draw=none, fill=none}]   {   E4}{  0}{ 11   }
\ganttbar[bar label font=\tiny\color{black},]{3}{5}{7}\\
\ganttgroup[inline=false,group label font=\scriptsize, group/.style={draw=none, fill=none}]   {   E5}{  0}{ 11   }
\ganttbar[bar label font=\tiny\color{black},]{2}{5}{6}\\
\ganttgroup[inline=false,group label font=\scriptsize, group/.style={draw=none, fill=none}]   {   E6}{  0}{ 11   }
\ganttbar[bar label font=\tiny\color{black},]{5}{5}{9}\\
\ganttvrule[
vrule/.append style={blue}]{$H = 10$}{9}
\end{ganttchart}
\end{minipage}

\end{figure}

\subsection{Reduction}\label{reduct}

 The  decision version of the {\detprobname} is as follows:\\
\textcolor{white}{..}\\
\fbox{\begin{minipage}{\dimexpr\textwidth-2\fboxsep-2\fboxrule\relax}
\textbf{INSTANCE }: { Finite vector $\mathbf{r}$ with duration $r_i \in Z^+$, for each $i\in \mathcal{E} = \{1,..,M\}$, a target horizon $H$  and a target number of potential bumps $B^*$ . }\\
\textbf{QUESTION }: { Does a feasible  schedule $S$ exist with at most $B^*$ potential bumps for a given instance?}  
\end{minipage}}
\\
\textcolor{white}{..}\\
To prove that the {\detprobname} is  $\mathcal{NP}-complete$, we reduce from the  \textsc{Subset-Sum} problem.

The  \textsc{Subset-Sum} problem is defined as follows using (\citet{karp1972reducibility}, \citet{garey1979computers}).\\
\textcolor{white}{aa}\\
\fbox{\begin{minipage}{\dimexpr\textwidth-2\fboxsep-2\fboxrule\relax}
\textbf{INSTANCE }: { Finite set A, size $a_j \in Z^+$, for each $j\in A$ and a target positive integer $W $. }\\
\textbf{QUESTION }: { Is there a subset $A^{*} \subseteq A$ such that the sum of sizes in $A^{*}$ is exactly $W$ i.e. $\sum_{j \in A^{*}} a_j = W$?}  
\end{minipage}}
\\
\textcolor{white}{..}\\
Given $A = \{1,\dots, N\}$, the sizes $[a_j]_{j \in A}$, and a target  $W$, let $\ntpss$ denote the reduced instance.   Let the total number of employees be  $M  =  N + \sum_{j \in A}a_j + 1$. The letter $i$ is used to denote the seniority level of an employee. We define different sets of employees and their seniority as described in \Cref{sets_emps}. These three sets represent the \textbf{critical}, \textbf{stable}, and the specific \textbf{last} employee. A critical employee is the one who will only cause bumps. A stable set employee will only suffer a bump. There is one critical employee $i_k$ for each index $k \in A$. Also, any two critical employees $i_k, i_{k+1}$ are separated by a block of stable set employees $\mathcal{E}^S_{k} $ of size $a_k$. There are $|A|$ blocks of stable employees.  Finally, we have a singleton set for the most junior employee to complete.
\begin{table}[!ht]
    \centering
    \caption{Employee sets for the reduction}
\label{sets_emps}
\parbox{1\linewidth}{
\centering
\resizebox{\linewidth}{!}{
  {\def\arraystretch{0.8}  
\begin{tabular}{llll}
\hline
\textbf{Sets}& & &\\
\hline
$\mathcal{E}$ & & set of all employees in order of seniority, = $\{1,...M\}$ & $|\mathcal{E}| = M$\\
$\mathcal{E}^C $ & &  set of seniority indices of critical employees,  = $\{i_k: i_k = k + \sum_{j = 1}^{k - 1}a_j, k\in A \}$ &$|\mathcal{E}^C| = N$ \\
$\mathcal{E}^S_{k} $ & &  $k^{th}$  block of stable employees,  = $\{i: i_k < i \leq i_k + a_k\}$ & $|\mathcal{E}^S_{k}| =  a_k$ \\
$\mathcal{E}^S $ & &  set of seniority indices of stable employees,  = $\cup_{k\in A} \mathcal{E}^S_{k}$ & $|\mathcal{E}^S| =  \sum_{k \in A}a_k$ \\
$\mathcal{E}^M $ & &  singleton set for the last employee,  = $\{M\}$ &$|\mathcal{E}^M| = 1$\\
\hline

\end{tabular}}}}

\end{table}

 The response delay for each employee set is given by  \Cref{rdtoemp}.  For any employee $i_k$ in the critical set, the response delay is simply the sum of all sizes till index $k$. For the stable employees in the $k^{th}$ block, the response delay is the same as the critical employee $i_{k-1}$. For the last employee,  it is simply the sum of all sizes.
\begin{equation}
\label{rdtoemp}
    r_{i} = \begin{cases}
  \sum\limits_{j = 1}^{k}a_j  & \text{ if } i = i_k \in \mathcal{E}^C, k \in A,  \\
  r_{i_{k-1}} & \text{ if } i \in \mathcal{E}^S,  i_k<i < i_{k+1} , k \in A, \\
  \sum\limits_{j \in A}a_j & \text{ if } i \in \mathcal{E}^M;
\end{cases}
\end{equation}
where $r_{i_0} = 0$.

The total planning horizon is set to $ H = C^*_0 - W$ where $C^*_0 = 2 \sum_{j \in A}a_j$ using \Cref{prop_makespan_r}. 
It is assumed that  $W \leq \sum_{j \in A} a_j$,  else it is polynomial to check that the  \textsc{Subset-Sum} is infeasible when $W > \sum_{j \in A} a_j$.
One can easily verify that the above instance has size $r_i \in \mathbb{N}$, for each $i \in \mathcal{E}$ and a target positive integer $H  \in \mathbb{N}$  and is achieved in polynomial operations. The following example illustrates the reduction.

\paragraph{\textbf{Example}:} Consider a  \textsc{Subset-Sum} instance with  $A = \{1,2,3\} $, sizes $a_1 = 1, a_2 = 4, a_3 = 7$ and  with $W = 5$. As a consequence $|\mathcal{E}^C| = 3, |\mathcal{E}^S| = 12, |\mathcal{E}^M| = 1, M = 16$. The seniority indexes in the critical,  stable, and last employee sets are defined below. 
\begin{itemize}
    \item $\mathcal{E}^C = \{1,3,8\}$
    \item $\mathcal{E}^S = \{2,4,5,6,7,9,10,11,12,13,14,15\}$
    \item $\mathcal{E}^M = \{16\}$
\end{itemize}
The response delay vector $R = [1,0,5,1,1,1,1,12,5,5,5,5,5,5,5,12]$. The horizon $H  = C^*_0 - W = 19$, where $C^*_0 = 24$. Color codes are used to denote employees from their respective set on the $x$-axis  using the following convention: \textcolor{magenta}{$\mathcal{E}^C $}, \textcolor{teal}{$\mathcal{E}^S $}, \textcolor{orange}{$\mathcal{E}^M $ }.  The NBS and the optimal schedule $S^*$ are in \Cref{gantt2} (a) and (b). 

\begin{figure}[ht!] 
\caption{(a) \textbf{NBS} : $C^{*}_{0} = 24$ \hspace{3cm} (b) $\mathbf{S^*}$ : $\mathbf{B(S^*) = 5}$ and $C(S^*) = 19$;  }
\label{gantt2}
\begin{minipage}{0.5\textwidth}

\begin{ganttchart}[
x unit  = 0.3cm,
hgrid = true,
vgrid=dashed, bar/.append style={fill=gray!50},
y unit chart=10,
vrule/.style={very thick, blue},
title height=.4,
title height=1,
y unit title=0.4cm,
bar height=0.6,
bar top shift = 0.2,
inline,
link/.style={-latex, draw=red, fill=red}
]{0}{24}

\gantttitle{\textbf{Notification Schedule }}{25} \\
\gantttitlelist[title label font=\tiny,]{1,...,25}{1}\\
\ganttgroup[inline=false,group label font=\scriptsize, group/.style={draw=none, fill=none}]   {  \textcolor{magenta}{ E1}}{  0}{ 11   }
\ganttbar[bar label font=\tiny\color{black},]{1}{0}{0} \\
\ganttgroup[inline=false,group label font=\scriptsize, group/.style={draw=none, fill=none}]   {  \textcolor{teal}{E2}}{  0}{ 11   }
\ganttbar[bar label font=\tiny\color{black},]{0}{2}{-1}\\
\ganttgroup[inline=false,group label font=\scriptsize, group/.style={draw=none, fill=none}]   {   \textcolor{magenta}{E3}}{  0}{ 11   }
\ganttbar[bar label font=\tiny\color{black},]{5}{1}{5}\\
\ganttgroup[inline=false,group label font=\scriptsize, group/.style={draw=none, fill=none}]   {   \textcolor{teal}{E4}}{  0}{ 11   }
\ganttbar[bar label font=\tiny\color{black},]{1}{5}{5}\\
\ganttgroup[inline=false,group label font=\scriptsize, group/.style={draw=none, fill=none}]   {   \textcolor{teal}{E5}}{  0}{ 11   }
\ganttbar[bar label font=\tiny\color{black},]{1}{5}{5}\\
\ganttgroup[inline=false,group label font=\scriptsize, group/.style={draw=none, fill=none}]   {   \textcolor{teal}{E6}}{  0}{ 11   }
\ganttbar[bar label font=\tiny\color{black},]{1}{5}{5}\\
\ganttgroup[inline=false,group label font=\scriptsize, group/.style={draw=none, fill=none}]   {  \textcolor{teal}{E7}}{  0}{ 11   }
\ganttbar[bar label font=\tiny\color{black},]{1}{5 }{5} \\
\ganttgroup[inline=false,group label font=\scriptsize, group/.style={draw=none, fill=none}]   {   \textcolor{magenta}{E8}}{  0}{ 11   }
\ganttbar[bar label font=\tiny\color{black},]{12}{5}{16}\\
\ganttgroup[inline=false,group label font=\scriptsize, group/.style={draw=none, fill=none}]   {   \textcolor{teal}{E9}}{  0}{ 11   }
\ganttbar[bar label font=\tiny\color{black},]{5}{12}{16}\\
\ganttgroup[inline=false,group label font=\scriptsize, group/.style={draw=none, fill=none}]   {   \textcolor{teal}{E10}}{  0}{ 11   }
\ganttbar[bar label font=\tiny\color{black},]{5}{12}{16}\\
\ganttgroup[inline=false,group label font=\scriptsize, group/.style={draw=none, fill=none}]   {   \textcolor{teal}{E11}}{  0}{ 11   }
\ganttbar[bar label font=\tiny\color{black},]{5}{12}{16}\\
\ganttgroup[inline=false,group label font=\scriptsize, group/.style={draw=none, fill=none}]   {  \textcolor{teal}{E12}}{  0}{ 11   }
\ganttbar[bar label font=\tiny\color{black},]{5}{12}{16}\\
\ganttgroup[inline=false,group label font=\scriptsize, group/.style={draw=none, fill=none}]   {   \textcolor{teal}{E13}}{  0}{ 11   }
\ganttbar[bar label font=\tiny\color{black},]{5}{12}{16} \\
\ganttgroup[inline=false,group label font=\scriptsize, group/.style={draw=none, fill=none}]   {   \textcolor{teal}{E14}}{  0}{ 11   }
\ganttbar[bar label font=\tiny\color{black},]{5}{12}{16}\\
\ganttgroup[inline=false,group label font=\scriptsize, group/.style={draw=none, fill=none}]   {   \textcolor{teal}{E15}}{  0}{ 11   }
\ganttbar[bar label font=\tiny\color{black},]{5}{12}{16}\\
\ganttgroup[inline=false,group label font=\scriptsize, group/.style={draw=none, fill=none}]   {   \textcolor{orange}{E16}}{  0}{ 11   }
\ganttbar[bar label font=\tiny\color{black},]{12}{12}{23}\\
\ganttvrule[
vrule/.append style={red}]{$H = 19$}{18}
\end{ganttchart}
\end{minipage}\hfill
\begin{minipage}{0.5\textwidth}
\begin{ganttchart}[
x unit  = 0.3cm,
hgrid = true,
vgrid=dashed, bar/.append style={fill=gray!50},
y unit chart=10,
vrule/.style={very thick, blue},
title height=.4,
title height=1,
y unit title=0.4cm,
bar height=0.6,
bar top shift = 0.2,
inline,
link/.style={-latex, draw=red, fill=red}
]{0}{24}

\gantttitle{\textbf{Notification Schedule}}{25} \\
\gantttitlelist[title label font=\tiny,]{1,...,25}{1}\\
\ganttgroup[inline=false,group label font=\scriptsize, group/.style={draw=none, fill=none}]   {  \textcolor{magenta}{ E1}}{  0}{ 11   }
\ganttbar[bar label font=\tiny\color{black},]{1}{0}{0} \\
\ganttgroup[inline=false,group label font=\scriptsize, group/.style={draw=none, fill=none}]   {  \textcolor{teal}{E2}}{  0}{ 11   }
\ganttbar[bar label font=\tiny\color{black},]{0}{0}{-1}\\
\ganttgroup[inline=false,group label font=\scriptsize, group/.style={draw=none, fill=none}]   {   \textcolor{magenta}{E3}}{  0}{ 11   }
\ganttbar[bar label font=\tiny\color{black},]{5}{0}{4}\\
\ganttgroup[inline=false,group label font=\scriptsize, group/.style={draw=none, fill=none}]   {   \textcolor{teal}{E4}}{  0}{ 11   }
\ganttbar[bar label font=\tiny\color{black},,bar/.append style={ fill=red}]{1}{0}{0}\\
\ganttgroup[inline=false,group label font=\scriptsize, group/.style={draw=none, fill=none}]   {   \textcolor{teal}{E5}}{  0}{ 11   }
\ganttbar[bar label font=\tiny\color{black},,bar/.append style={ fill=red}]{1}{0}{0}\\
\ganttgroup[inline=false,group label font=\scriptsize, group/.style={draw=none, fill=none}]   {   \textcolor{teal}{E6}}{  0}{ 11   }
\ganttbar[bar label font=\tiny\color{black},,bar/.append style={ fill=red}]{1}{0}{0}\\
\ganttgroup[inline=false,group label font=\scriptsize, group/.style={draw=none, fill=none}]   {  \textcolor{teal}{E7}}{  0}{ 11   }
\ganttbar[bar label font=\tiny\color{black},,bar/.append style={ fill=red}]{1}{0}{0} \\
\ganttgroup[inline=false,group label font=\scriptsize, group/.style={draw=none, fill=none}]   {   \textcolor{magenta}{E8}}{  0}{ 11   }
\ganttbar[bar label font=\tiny\color{black},]{12}{0}{11}\\
\ganttgroup[inline=false,group label font=\scriptsize, group/.style={draw=none, fill=none}]   {   \textcolor{teal}{E9}}{  0}{ 11   }
\ganttbar[bar label font=\tiny\color{black},]{5}{7}{11}\\
\ganttgroup[inline=false,group label font=\scriptsize, group/.style={draw=none, fill=none}]   {   \textcolor{teal}{E10}}{  0}{ 11   }
\ganttbar[bar label font=\tiny\color{black},]{5}{7}{11}\\
\ganttgroup[inline=false,group label font=\scriptsize, group/.style={draw=none, fill=none}]   {   \textcolor{teal}{E11}}{  0}{ 11   }
\ganttbar[bar label font=\tiny\color{black},]{5}{7}{11}\\
\ganttgroup[inline=false,group label font=\scriptsize, group/.style={draw=none, fill=none}]   {  \textcolor{teal}{E12}}{  0}{ 11   }
\ganttbar[bar label font=\tiny\color{black},]{5}{7}{11}\\
\ganttgroup[inline=false,group label font=\scriptsize, group/.style={draw=none, fill=none}]   {   \textcolor{teal}{E13}}{  0}{ 11   }
\ganttbar[bar label font=\tiny\color{black},]{5}{7}{11} \\
\ganttgroup[inline=false,group label font=\scriptsize, group/.style={draw=none, fill=none}]   {   \textcolor{teal}{E14}}{  0}{ 11   }
\ganttbar[bar label font=\tiny\color{black},]{5}{7}{11}\\
\ganttgroup[inline=false,group label font=\scriptsize, group/.style={draw=none, fill=none}]   {   \textcolor{teal}{E15}}{  0}{ 11   }
\ganttbar[bar label font=\tiny\color{black},]{5}{7}{11}\\
\ganttgroup[inline=false,group label font=\scriptsize, group/.style={draw=none, fill=none}]   {   \textcolor{orange}{E16}}{  0}{ 11   }
\ganttbar[bar label font=\tiny\color{black},]{12}{7}{18}\\
\ganttvrule[
vrule/.append style={blue}]{$H = 19$}{18}
\end{ganttchart}

\end{minipage}\hfill
\end{figure}

Let us define a type of Notification schedule that will help us establish the proof.
\begin{definition}
Block Schedule $S(A^\prime)$:  Given $A^\prime \subseteq A$, a block schedule $S(A^\prime) = [(s_i, e_i)]_{i\in \mathcal{E}}$ for $\ntpss$   is constructed as follows:
\begin{equation}
\label{sstontp}
     s_{i} = \begin{cases}
  0  & \text{if } i = 1,  \\
  s_{i_k} + r_{i_k} - r_i & i = i_k + 1, \forall k \in A \backslash A^\prime, \\
  s_{i-1}& \text{else } ;
\end{cases}
 \end{equation}
\end{definition}
Alternately, all employees of the stable set correspond to $k \in A^\prime$.
One can easily verify that $\forall k \in A, r_{i_k} - r_i = r_{i_k} - r_{i_{k-1}} = a_k $ when $i = i_k + 1$ in \Cref{sstontp}.  \Cref{gantt2} (b) is an example of the block schedule~$S(A^\prime)$ with  $A^\prime = \{1,2\}$.

\subsection{Intuition}\label{intuit}
The crux of the reduction lies in how the response delays are defined. They ensure the following.

\begin{itemize}
    \item There are no bumps between any two employees in the stable set because for any $i,i^\prime \in \mathcal{E}^S$ such that $i < i^\prime$ we have $r_i \leq r_{i^\prime}$.
    \item There are no bumps between any two critical set employees, as we have $r_{i_k} \leq r_{i_{k^\prime}}$ for any $i_k, i_{k^\prime} \in \mathcal{E}^C$ such that $i_k < i_{k^\prime}$.
    \item Any critical employee $i_k$ can potentially bump only a whole block of stable set employees $\mathcal{E}^S_k$ to create time savings or bump none. This means the critical employee $i_k$ can introduce a total of $a_k$ potential bumps in the schedule. Consequently, the total potential bumps $B(S)$ will be a sum of these sizes. The aim is now to introduce bumps so that  $B(S)$ equals $W$. For example, in \Cref{gantt2} (a), one can see that E3 corresponding to index 2, can only bump E4, E5, E6, and E7 corresponding to the size $a_2 = 4$. E3 cannot bump any employee junior to E7.  It is also ensured that no stable employee bumps any other employee from the entire set. 
    \item A block schedule is always feasible and optimal. The block schedule $S(A)$, which notifies everyone at the start, is always feasible.  Furthermore, when the total potential bumps $B(S) = \sum\limits_{k\in A^\prime} b_{i_k}$, where $A^\prime \subset A$, the block schedule~$S(A^\prime)$ gives the minimum makespan and hence is also optimal.
    \item Any schedule $S$ such that $B(S) < W$ is infeasible. If a critical employee $i_k$ potentially bumps a block of stable set employee $ \mathcal{E}^S_k$, then it also creates maximum time savings of $a_k$. Note that a block schedule will create time savings of exactly $a_k$ units for a potential bump by critical employee $i_k$. Thus, starting from the NBS schedule with a makespan of $C^*_0$, we set $H = C^*_0 - W$ to ensure that any feasible solution has at least $W$ potential bumps.

\end{itemize}

The goal is to find a feasible block schedule $S$, such that $B(S) = W$. When this happens, all the indexes $\{k_1, k_2,\dots, k_Z\}$ such that critical employee $ I_{k_u}$, causes all the potential bumps, $\forall u \in \{1,\dots,Z\}$  are part of the subset, and we would have a YES answer to the  \textsc{Subset-Sum} problem.

\subsection{Proof} \label{proof}
Our main result in this section establishes the computational complexity of the {\detprobname}. The proofs of all propositions, corollaries, lemmas, and theorems are given in the attached E-Companion.
Let $\feasS\bigg([a_j]_{j \in \mathcal{A}}, W\bigg)$ and $\optS\bigg([a_j]_{j \in \mathcal{A}}, W\bigg)$ denote the feasible and optimal set of solutions for $\ntpss$. 

\begin{proposition}
 The block schedule $S(A^\prime)$ has the following properties:
\label{proper:1}
\begin{enumerate}
    \item $\forall k \in A, s_{i_k} = s_{i_k - 1}$.
    \item $\forall i \in \mathcal{E}^S_{k}, \forall k \in A$,  $e_i = e_{i_k}$ or $s_i = s_{i_k}$.
    \item $s_{i_{k+1}} = s_{i_{k}} + \mathbf{1}_{k \notin A^\prime}a_{k}$, where $\mathbf{1}$ denotes the indicator variable.
    \item $B(S(A^\prime)) = \sum\limits_{j \in A^\prime} a_j$
    
    \item $C(S(A^\prime)) = C^*_0 - \sum\limits_{j \in A^\prime } a_j $.
\end{enumerate}
\end{proposition}
Properties 1 and 2 mean that the block schedule respects the seniority constraint of the problem. Given the notification time of any critical employee $i_k$, property 3 allows us to find the notification time of the next critical employee $i_{k+1}$. Property 4 gives an equation to compute the total potential bumps suffered by the block schedule.  Property 5 establishes the total makespan of the block schedule.  

\begin{proposition}
\label{prop:1}
For any $\ntpss$, $C = e_M$.
\end{proposition}

\begin{proposition}
\label{prop:2}
$\forall k \in A$, $\forall l \in A$ employee $i_k$ cannot bump employee $i_l$.

\end{proposition}

\begin{proposition}
\label{prop:3}
 $\forall i \in \mathcal{E^S}$, $\forall i^\prime \in \mathcal{E}$ employee $i$ cannot bump employee $i^\prime$.
\end{proposition}

\begin{proposition}
\label{prop:4}
$\forall k \in A $, $\forall i  \in \mathcal{E}: i \geq i_{k+1}$, employee $i_k$ cannot bump employee $i$.
\end{proposition}

\begin{corollary}
\label{cor:1}
$\forall k \in A$, $\forall i\in \mathcal{E}$ if employee $i_k$  bumps employee $i$, then $i \in \mathcal{E}^S_{k}$.
\end{corollary}

 \Cref{prop:1} defines the makespan of the schedule. Since the last employee $i_M$ has the maximum response delay, it is simple to show that the end of its interaction, $e_M$,  marks the end of the schedule for $\ntpss$. \Cref{prop:2}, \Cref{prop:3}, and \Cref{prop:4} establish that it is only the critical set employees who cause the bumps. Moreover,  any critical employee $i_k$ can only bump the immediately following $a_k$ employees in the order of seniority, all of whom belong to the stable set from \Cref{cor:1}.

\begin{lemma}
\label{lem:1}
 $\forall S^* \in \optS\big([a_j]_{j \in \mathcal{A}}, W\big) $, $\forall k \in A, b_{i_k} = a_k $ or $0$. 
\end{lemma}

 \begin{corollary}
\label{cor:2}
 $\forall S \in \feasS\big([a_j]_{j \in \mathcal{A}}, W\big)$, the total bumps $B(S) = \sum\limits_{k\in A} b_{i_k}$.
\end{corollary}

\begin{lemma}
\label{lem:2}
 $S(A) \in \feasS\big([a_j]_{j \in \mathcal{A}}, W\big)$. 
\end{lemma}
 \Cref{lem:1} proves that any critical employee $i_k$ will bump the immediately following $a_k$ junior employees or bump no one.  As a consequence of \Cref{lem:1}, \Cref{cor:2} establishes the total bumps in the schedule. Next, we show that the block schedule $S(A)$ is always feasible for the $\ntpss$ in \Cref{lem:2}. Hence, the reduced problem is always feasible.

\begin{lemma}
\label{lem:3}
$\exists A^\prime \subseteq A$ such that $S(A^\prime) \in \optS\big([a_j]_{j \in \mathcal{A}}, W\big)$. 
\end{lemma}
\begin{corollary}
\label{cor:3}
 $\forall S^* \in \optS\big([a_j]_{j \in \mathcal{A}}, W\big)$, the optimal objective value $B(S^*) \geq W$.
\end{corollary}

The final piece of the proof comes from \Cref{lem:3}. This means that given any horizon $W$, some subset $A^\prime$ and an optimal block schedule $S(A^\prime)$ exist. This is a significant result as it allows only for concentrating on a solution space that comprises block schedules.
\begin{theorem}
\label{theo}
For $\ntpss$, $ B(S^*) = W \Leftrightarrow \exists A^* \subseteq A $ such that $\sum_{j \in A^*} a_j = W$.
\end{theorem}

\begin{theorem}
\label{theo2}
  \textsc{NTP}$(I)$ is \emph{NP-complete}.
\end{theorem}

 \Cref{theo} proves that for any $\ntpss$, the number of bumps in the optimal schedule $B(S^*)$ is always equal to the target $W$ defined by the  \textsc{Subset-Sum} if and only if the  \textsc{Subset-Sum} problem is feasible. Hence, the offline problem is  $\mathcal{NP}$-complete from  \Cref{theo2}, which marks the most important contribution of the paper. The hardness of the problem stems from the fact that one needs to count how many bumps an employee can cause. 
 
 \subsection{Comparing offline and online solutions}
Consider the online version of the problem, where response delays are unknown but all other assumptions of {\detprobname} hold. Suppose we have an online algorithm that observes the system at each discrete-time epoch and determines when to notify the next employee. Naturally, we want to evaluate the performance of this algorithm. Competitive ratios (\cite{borodin2005online}) are a useful measure here, as they compare the quality of an online algorithm’s solution to the optimal result of the offline model, where full information is available. The competitive ratio is simply the worst-case ratio between the cost of the algorithm's solution and the cost of an optimal offline solution. The following theorem applies.

 \begin{theorem}
\label{theo4}
Given any online algorithm, there exists an instance $I$ in which the algorithm suffers either the maximum potential bumps or has vacant shifts. The offline counterpart will have zero potential bumps and zero vacancies.
\end{theorem}

In this case, the competitive ratio is undefined. The proof is attached in E-Companion. As a consequence of \Cref{theo4}, the performance of an online algorithm can be arbitrarily poor. This is a significant result as it suggests the problem remains hard on-line.

%% file: dntp.tex
\section{Dynamic Notification Timing Problem}\label{dntp}
In reality, employee response delays are stochastic and can vary widely. However, we assume the distribution of these delays is known. This uncertainty adds complexity, as there’s a risk that some employees may not respond at all—specifically if their response delay exceeds the planning horizon. In such cases, waiting for responses to avoid bumps might result in unfilled shifts. We define the dynamic stochastic version of {\detprobname}, referred to as {\stochprobname}, in \Cref{desc_dntp} and present a two-stage stochastic formulation to compute an offline policy. In \Cref{PFA}, we develop approximate policy functions for {\stochprobname}.

\subsection{Problem Formulation - {\stochprobname}} \label{desc_dntp}

All data parameters of the deterministic problem, {\detprobname}, apply to its stochastic dynamic counterpart, {\stochprobname}. This problem represents the actual online problem faced. A shorter bump cutoff time, $D$, is introduced to allow senior employees to retain their advantage to bump junior ones. To account for uncertainty, the employee pool exceeds the number of shifts, creating a buffer if some employees do not respond in time. Employee response delays, $r_i$, are modeled as independent, identically distributed random variables. As in the deterministic version, all employees are eligible for any shift, though a cap restricts simultaneous notifications for operational efficiency.

 The {\stochprobname} is a multi-stage optimization problem and is really difficult to solve. We propose a two-stage stochastic model to solve it: in the first stage, we determine the notification schedule for each employee; in the second stage, bumps are counted based on the realization of employee response. Let ${\Omega}$ denote a set of scenarios and $\omega \in \Omega$ to denote a scenario with the corresponding response delay as $\mathbf{r}_\omega$. The MIP formulation for this problem is given by ${\textbf{\textsc{MIP}}_{\textbf{DNTP-S}}} \eqref{objmain} - \eqref{eqws13}$. Table \ref{notation2} defines the additional variables and parameters.

\begin{table}[ht]

\parbox{1\linewidth}{
\centering
\caption{ Parameters and Variables of the \textsc{DNTP}}
\label{notation2}}
\resizebox{\linewidth}{!}{
   {\def\arraystretch{1.2} 
\begin{tabular}{llll}
\hline
\textbf{Parameters}&&\\
\hline
$D$ & cutoff time to bump a junior employee\\
$G$ & reward for filling a shift\\
$W$ & maximum number of concurrent notifications \\
\hline
\textbf{Variables}&&\\
\hline
${y}_{ij\omega}\in \{0,1\}$ &  1 if employee $i \in \mathcal{E}$ bumps employee $j \in \mathcal{E}$ in scenario $\omega$ else 0\\
$\bar{z}_{i\omega}\in \{0,1\}$ & 1 if employee $i \in \mathcal{E}$ responds within the horizon and gets a shift in scenario $\omega$ else 0\\
$\hat{z}_{i\omega}\in \{0,1\}$ &  1 if employee $i \in \mathcal{E}$ responds within the horizon and does not get a shift in scenario $\omega$ else 0\\
$\theta_\omega\in Z^+ $ & number of shifts vacant in scenario $\omega$\\ 
\hline
\end{tabular}}}
\end{table}
\resizebox{16cm}{!}{
\centering
\makebox[\textwidth]{\begin{minipage}{\dimexpr\textwidth-4\fboxsep-1\fboxrule\relax}
\begin{equation}
{\textbf{\textsc{MIP}}_{\textbf{DNTP-S}}}:= \quad  \min \quad   \frac{1}{|\Omega|}\sum\limits_{\omega \in \Omega} Q(s, \mathbf{r}_\omega) 
\label{objmain}
\end{equation}
\begin{align}
\label{eqws2} \mbox{s.t. }\quad \quad & s_i \leq s_j  & \forall i,j\in \mathcal{E}, i < j \\
\label{eqws9}&s_i + 1= s_{i + W} & \forall i, i + W \in \mathcal{E}\\
\label{eqws14}& s_{i} \in {Z}^+ & \forall i\in \mathcal{E}\\
\end{align}

\begin{equation}
\text{where } \quad Q(s, \mathbf{r}_\omega):= \quad  \min \quad    \bigg\{ G\theta_\omega  + \sum_{i\in \mathcal{E}}\sum_{j\in \mathcal{E}:i<j} y_{ij\omega} \bigg\} \\
\end{equation}
\begin{align}
 \mbox{s.t. }\quad  &s_i + r_{i\omega} \geq  (H+1)(1 - \bar{z}_{i\omega} - \hat{z}_{i\omega}) & \forall i\in \mathcal{E} \label{eqws3} \\
& s_i + r_{i\omega} \leq  H + r_{i\omega} (1 - \bar{z}_{i\omega} - \hat{z}_{i\omega}) & \forall i\in \mathcal{E} \label{eqws4}\\
& s_i - s_j + \delta_{ij\omega} \leq \delta_{ij\omega}y_{ij\omega} +  (H + r_{i\omega})(1 - \bar{z}_{i\omega} + \hat{z}_{i\omega}) & \forall i,j\in \mathcal{E}, i<j,  r_{j\omega} \leq r_{i\omega} \leq D   \label{eqws5}\\
&\sum_{i\in \mathcal{E}}\bar{z}_{i\omega} + \theta_\omega \geq L &\label{eqws6}\\
&\sum_{j<i}\hat{z}_{j\omega}  \leq i( 1 - \bar{z}_{i\omega}) &\forall i \in \mathcal{E} \label{eqws10}\\
&y_{ij\omega} \in \{0,1\} & \forall i,j\in \mathcal{E}, i < j\label{eqws7}\\
&\bar{z}_{i\omega}, \hat{z}_{i\omega} \in \{0,1\} & \forall i\in \mathcal{E}\label{eqws12}\\
&\theta_\omega \geq 0 & \label{eqws13} 
\end{align}
\end{minipage}}}
\vspace{0.5cm}\\

We use a weighted objective function that aims to minimize both bumps and shift vacancies. That is, instead of enforcing a strict shift-vacancy constraint as in ${\textbf{\textsc{MIP}}_{\textbf{NTP}}}$, we introduce a high penalty, $G$, for any unfilled shifts.    In the offline version of this problem ({\detprobname}), it is straightforward to determine the number of unfilled shifts. Any employee with a response time shorter than the planning horizon can fill a shift, meaning vacant shifts occur only when fewer employees respond within the horizon than available shifts. Thus, the notification schedule can be optimized to minimize bumps while satisfying the shift-vacancy constraint. However, in the stochastic case, since we need a unified notification schedule that works across all possible scenarios, the number of vacant shifts cannot be predetermined due to varying response times across scenarios. Given the first-stage decisions \( s_i \), the number of employees responding back will also vary in each scenario. In some scenarios, all employees may need to be notified to fill the shifts, while in others, all shifts may be scheduled before notifying every employee. Thus, each scenario will have a different number of vacancies and bumps.

To correctly count the bumps, we introduce two new variables, \( \bar{z}_{i\omega} \) and \( \hat{z}_{i\omega} \), to track whether employee \( i \) has responded and whether they are assigned a shift, respectively for each realized scenario of $\mathbf{r}_\omega$. Constraint \eqref{eqws2} ensures that seniority is respected, while Constraint \eqref{eqws9} enforce the operational limit of notifying a maximum of \( W \) employees per epoch.  The second-stage constraints are provided in equations \eqref{eqws3} - \eqref{eqws12}. Constraint \eqref{eqws3} and \eqref{eqws4} determine when employee \( i \) responds. Constraint \eqref{eqws5} tracks potential bumps, and constraint \eqref{eqws6} counts vacant shifts in each scenario. Additionally, constraint \eqref{eqws10} ensures that if employee \( i \) responds and is assigned a shift (i.e., \( \bar{z}_{i\omega} = 1 \)), then all senior employees who responded earlier will also be assigned a shift (i.e., \( \hat{z}_{i\omega} = 0 \)). Constraints \eqref{eqws6}, \eqref{eqws10} together ensure that \( \bar{z}_{i\omega} \), \( \hat{z}_{i\omega} \) are assigned appropriate values if an employee gets a shift or not. We introduce a large constant \( G \) to ensure complete shift assignment. Note that \( G \) is used only to prioritize shift assignments and does not reflect the actual cost of leaving a shift vacant. Estimating the true cost of a vacant shift is challenging, as management may resort to alternative means to fill the shift as a last resort.

This formulation is still quite complex and computationally expensive to solve. We also provide a complete information offline formulation for a single scenario $\omega$, denoted as \( {\textbf{\textsc{MIP}}_{\textbf{NTP2}}} \), which is detailed in the E-Companion. This offline formulation is a simplified version of ${\textbf{\textsc{MIP}}_{\textbf{DNTP-S}}}$ requiring only one of the variables \( \bar{z}_{i\omega} \), \( \hat{z}_{i\omega} \) to track responses.


\subsection{Designing Policy Function Approximations} \label{PFA}

To tackle the complexity of the {\stochprobname}, we focus on designing approximate policies using Policy Function Approximations (PFA) as outlined by \citet{powell2021reinforcement}. PFA is a policy search approach that is especially valuable when a straightforward, easy-to-implement policy is desirable. These policies provide a direct mapping from states to actions using an analytical function, thereby eliminating the need to solve complex optimization problems in each decision stage. Examples of PFAs include lookup tables, linear decision rules, monotone threshold policies, and nonlinear models. The central challenge is to develop a policy structure that is well-suited to the specific problem context, balancing simplicity with performance.

For the {\stochprobname}, we approximate optimal policies by employing a monotone threshold-based structure. A policy follows a monotone threshold approach when it selects an action \( a_k \) only up to a defined threshold \( \Lambda_k \) at decision epoch \( k \). Within this framework, the action \( a_k \) is directly influenced by the residual \( (\Lambda_k - \lambda_k) \), where \( \lambda_k \) is the current value of a target feature or basis in the state \( X_k \), and \( \Lambda_k \) is the optimal threshold for that feature. A typical example of such a policy is the classic \( (s, S) \) inventory policy, where restocking occurs only when inventory drops below \( s \), and replenishment brings it up to \( S \). In this case, the total inventory level serves as the feature driving the threshold. Thus, our task centers on designing effective features and determining their associated optimal thresholds $\Lambda_k$. However, finding these optimal thresholds is challenging. We adopt an approximate method to estimate appropriate threshold values to address this problem, optimizing two distinct policies for this problem. This approach allows us to balance computational feasibility with policy performance efficiently.

\begin{itemize}
\item  \textbf{Notify $\eta$ and Wait $w$ epochs} (NAW$(\eta,w)$) - A straightforward two-parameter static policy that sends \( \eta \) notifications every \( w \) epochs. The parameters \( \eta \) and \( w \) can be adjusted for optimal performance. This state-independent approach does not consider any information from the state and consistently notifies individuals at fixed intervals. Since the parameters remain unchanged over time, we can utilize a simple brute-force search to determine the optimal estimates by minimizing the average cost across all simulation runs. Our industrial partner currently employs this policy.

\item \textbf{Offline Notification Policy} (ONP) - ONP is a dynamic policy in which the threshold parameters vary over time. This heuristic policy depends solely on the elapsed time since the initiation of the notification process rather than the current state of the system.  The target feature utilized is the cumulative number of notifications issued at each decision epoch. Next, we will discuss the algorithm employed to estimate these thresholds.
\end{itemize}


\Cref{alg:two} presents pseudo-code for the function that defines estimates for thresholds \( \Lambda_k \) for the ONP. Here, \( \Lambda_k \) represents the cumulative number of notifications that should have been sent by epoch \( k \). The instance parameters, the set of offline instances \( \Omega \), and an aggregator function \( q \) are provided as inputs. Algorithm \ref{alg:two} first generates and solves a single scenario \( \omega \in \Omega \) using  ${\textbf{\textsc{MIP}}_{\textbf{NTP2}}}$. The solution obtained is then used to compute the feature value \( \lambda_k^{\omega} \) (e.g. the number of people notified) for each decision epoch \( k \) of scenario \( \omega \). We calculate these features for all instances in the scenario set. The aggregator function \( q \) is then applied to these feature values from the exact offline solutions, generating approximate thresholds \( \bar{\Lambda}_k = q([\lambda_k^{\omega}]_{\omega \in \{1,\dots,|\Omega|\}})\) across all instances for each decision epoch \( k \). The function outputs estimates of the threshold values for each decision epoch, represented as \( \bar{\Lambda} = [\bar{\Lambda}_{k}]_{k \in \mathcal{K}} \).
\begin{algorithm}
   \SetKwFunction{Function}{Compile}
   \SetKwProg{Fn}{Function}{:}{end}
\caption{An algorithm to estimate the threshold for all features}
\label{alg:two}

\Fn{\Function{$M, L, H, D, {\Omega}, q$}}{  
    \For{$\omega \in \{1,\dots,|{\Omega}|\}$}{
        Solve instance with realization $\mathbf{r}_\omega$ using ${\textbf{\textsc{MIP}}_{\textbf{NTP2}}}$\;
        Compute the value of each feature $\lambda^{\omega} = [\lambda_k^{\omega}]_{k \in \mathcal{K}}$;
    }
    
    \For{$k \in \mathcal{K}$}{
        $\bar{\Lambda}_k = q([\lambda_k^{\omega}]_{\omega \in \{1,\dots,|\Omega|\}})$\;
    }
    
    \Return{$\bar{\Lambda}$}
}
\end{algorithm}

In any decision epoch $k$, the decision $a = \lfloor\bar{\Lambda}_{k} - \lambda_k^\omega\rceil$ is the residual between the calculated target number of notifications $\bar{\Lambda}_{k}$ to be sent from offline solutions and its value in the current state, $\lambda_k^\omega$ for scenario $\omega$.   $\lfloor.\rceil$ denotes the nearest integer after rounding.

Algorithm \ref{alg:two} that designs the ONP offers a practical and computationally efficient heuristic solution to the stochastic problem ${\textbf{\textsc{MIP}}_{\textbf{DNTP-S}}}$.  Given a scenario set ${\Omega}$, one would use ${\textbf{\textsc{MIP}}_{\textbf{DNTP-S}}}$ to obtain a near-optimal policy. However, even for a small number of scenarios, this formulation explodes in size.  Algorithm \ref{alg:two} simplifies this complex formulation by decomposing the stochastic problem into smaller, manageable deterministic subproblems. Specifically, Algorithm \ref{alg:two} breaks down various stochastic scenarios into individual deterministic instances and leverages the relative simplicity of solving deterministic notification timing problems rather than directly addressing the full stochastic model. After deriving optimal or near-optimal solutions for each deterministic scenario, the results are aggregated to create a unified notification schedule that smooths out the uncertainty across the different instances considered similar to what the stochastic model will aim to do. Moreover, Algorithm \ref{alg:two} also allows the user to develop different threshold policies for different target features.

The primary goal of this approach is to develop a simple yet dynamic policy with time-varying threshold parameters. Various other threshold policies can be designed for different feature functions. In the literature, offline solutions have often been used to guide the design of effective online policies. For example, (\citet{pham2023prediction}) proposed a prediction-based approach for online dynamic radiotherapy scheduling, where a regression model is trained to map patient arrival patterns to optimal waiting times, leveraging offline solutions with full knowledge of future arrivals. Similarly, (\citet{de2021integrated}) tackled multi-stage optimization under uncertainty by integrating a two-stage offline strategy with an online greedy heuristic, achieving notable improvements in solution quality by enhancing offline/online integration. \citet{kong2022end} introduced a machine learning model that directly predicts the optimal action using offline solutions, using an energy-based parameterization of the model. For a comprehensive review of methods using offline full information solutions, we refer the reader to (\citet{sadana2024survey}, \cite{mandi2023decision}).

%% file: exp.tex
\section{Experiments}\label{results}
We conducted two sets of experiments to assess our approach. In the first set of experiments, we compare different policies with the weighted objective of potential bumps and vacant shifts. The second experiment aimed to evaluate our proposed methodology when preferences are different, representing the case of actual bumps. All experiments and models were implemented in Python, and we employed GUROBI 10.0 to solve the offline optimization problem instances. We set a maximum solving time of 4 minutes for each instance. GUROBI typically found the optimal solution quickly, often within seconds for all instances, though occasionally, it required more time to confirm optimality.

\subsection{Scenarios and Evaluation}
Experiments are conducted with a fixed planning horizon of 6 hours ($H = 6$ ) and a constant number of available shifts ($L = 50$). A group of 150 employees ($M = 150$) with identical preferences was available to fill these shifts. We examined two cutoff times: 2 hours and 3 hours (actual cutoff time). Employees are expected to be aware of the cutoff time and would respond before this time limit is reached.

Each experiment was defined by a tuple $(M, L, H, D)$, encompassing the number of employees, available shifts, the planning duration, and the cutoff time, respectively. This tuple was referred to as the setting for a specific experiment. We developed a simulation model to simulate the operation of the electronic system. This model accepts input parameters, including the notification policy, observes the current state of the system, and then takes action according to the input policy. Next, we define the number of instances used for each policy in the training and testing phase. We solve the actual stochastic formulation ${\textbf{\textsc{MIP}}_{\textbf{DNTP-S}}}$  for 200 scenarios. We attempted to solve the formulation for 500 scenarios but encountered memory issues. ONP is trained on a different set of 1000 instances and validated on 500 validation instances.  NAW is simulated on the 500 validation instances, and the best parameter is chosen. Finally, all policies are tested on the same 500 instances. Furthermore, a Notify-all (NA) policy was implemented, which notifies all employees at the beginning. This policy was included to emphasize its adverse impact on bumps.

We employ the following methodology to validate our dynamic policies:
\begin{itemize}
\item Threshold Estimation (Training): 1000 instances are used to estimate thresholds
\item Validation: Simulation of 500 validation instances to assess their performance.
\item Selection and Evaluation: Select the most effective aggregator functions based on their performance during validation. These selected functions are then evaluated using 500 test instances. 
\end{itemize}

This approach ensures that our chosen policy parameters can be generalized effectively when applied to previously unseen data. For the heuristic NAW, we explore all possible combinations of parameters, such as $\eta$ and $w$, using a predefined set of potential values. Each combination of parameters is simulated, and the best-performing parameter set is determined during this evaluation phase using the test set.

\subsection{Real World Data}
We evaluated our approach using real-world data from our industrial partner, analyzing shift scheduling data from 2018 to 2020. In some cases, the number of required shifts is known well in advance, allowing the on-demand scheduling system to begin up to four days before the shifts start. However, our experiments focus on a single day of operation with a six-hour planning window, specifically for shifts added the day before.

\begin{figure}[!ht]
    \centering

    \begin{minipage}{.5\textwidth}
        \centering
\caption{CDF plot of Response Delay}
\label{dist_rd}
\includegraphics[scale=0.54]{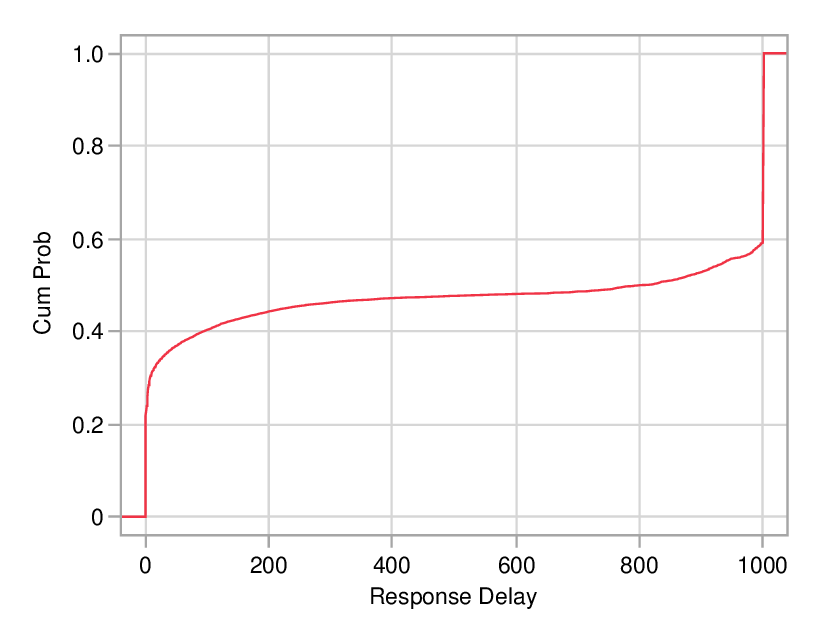}

     \end{minipage}%
    \begin{minipage}{.5\textwidth}
        \centering
    \caption{ Cumulative Notifications Sent}
    
    \label{detfigureall}
    \includegraphics[scale=0.2]{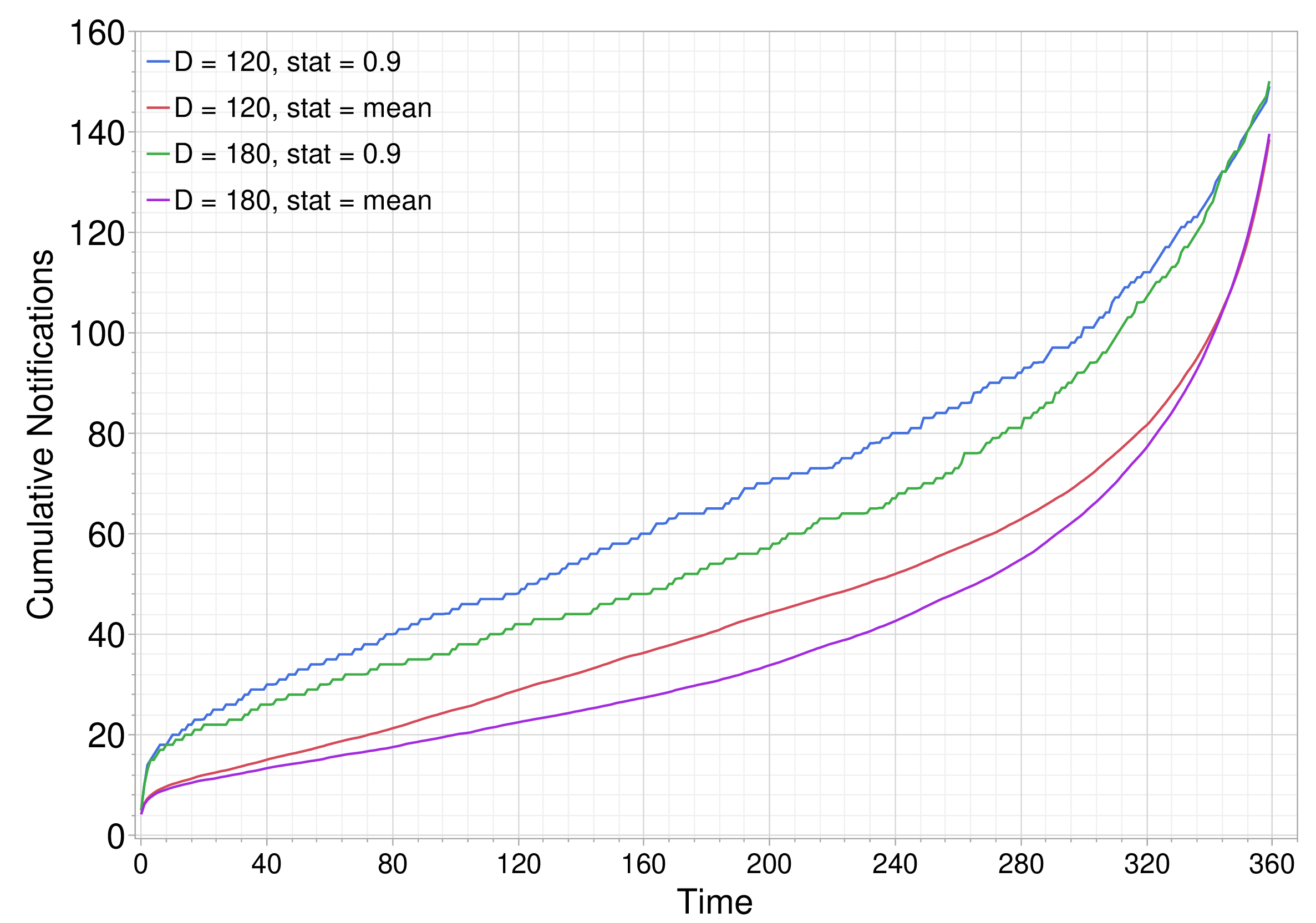}    \end{minipage}
    
    \vspace*{10pt}
\end{figure}

The current scheduling system uses a static NAW policy with fixed parameters, which we fine-tuned based on observed real-world response delays. After data cleaning, we analyzed 23,204 response delay records, as shown in. \Cref{dist_rd} presents the cumulative distribution of these response delays. Our analysis was limited to response times for shifts scheduled within a single day. We observed that only 50\% of employees responded to the system to select their shifts, with a large proportion doing so within one minute of receiving a system notification. This suggests those who want to work do indeed respond quickly as there is no incentive to wait and select shifts. It's important to note that the system allows employees to choose from remaining vacant shifts even after the planning period ends. However, the goal for management is to assign all shifts during the designated planning window. The data suggests that twice the number of employees is needed to fill all shifts, given that only 50\% respond. Interestingly, the response cutoff had little effect in our analysis. Response delays were sampled from the real-world data to generate instances split into training, validation, and testing sets. Additionally, the system limits notifications to 5 per epoch, so $W = 5$.

\subsection{Offline Solution Analysis}
 \Cref{alg:two} estimates thresholds of all features for 1000 training instances for each setting. We use descriptive statistics such as the ``mean" and ``percentiles" as aggregator functions for all features across all instances for each decision epoch.  \Cref{detfigureall} shows the variation of cumulative notifications concerning time when $D = 2, 3$ hours. Only the mean and the 90\% percentile were plotted for clarity. Notifications are sent near-linearly for both aggregator functions from the start. However, the rate changes as the horizon nears and notifications are sent faster. This is expected as employees notified in this period will not have sufficient time to bump others. In other words, employees who take a long time to respond are pushed out of the horizon to avoid bumps. The main difference between the two aggregators is that more notifications are sent earlier with 90\% aggregation. A shorter cutoff encourages more notifications to be sent with a similar profile.


 

\Cref{tab:offbumps} shows the average bumps and vacancy across all offline instances for the two cutoff values. It suggests that given prior knowledge of employee response times; one can almost eliminate bumps from the process without compromising on shifts being left vacant. This is mainly possible because many employees respond quickly or do not respond. One only needs to manage the few that take long to respond.
     \begin{table}
    
      \centering
      \caption{Average Optimal Bumps and Vacant Shifts in the offline solution}
\label{tab:offbumps}
    {\def\arraystretch{1.25}  
    \resizebox{4.5cm}{!}{
\begin{tabular}{|c|c|c|}
\hline
\textbf{D} & \textbf{Bumps} & \textbf{Vacant shifts} \\ \hline
\textbf{2} & 0.80          & 0                \\ \hline
\textbf{3} & 2.23           & 0                \\ \hline
\end{tabular}

}}
\end{table}

\subsection{Policy Evaluation - Same Preferences}

\subsubsection{Solving the stochastic formulation} 
\begin{figure}[!ht]
    \centering
\resizebox{7.5cm}{!}{
    \begin{minipage}{.48\textwidth}
        \centering
    \caption{Train vs Test performance of ${\textbf{\textsc{MIP}}_{\textbf{DNTP-S}}}$ - Average Bumps}
\includegraphics[width=\textwidth]{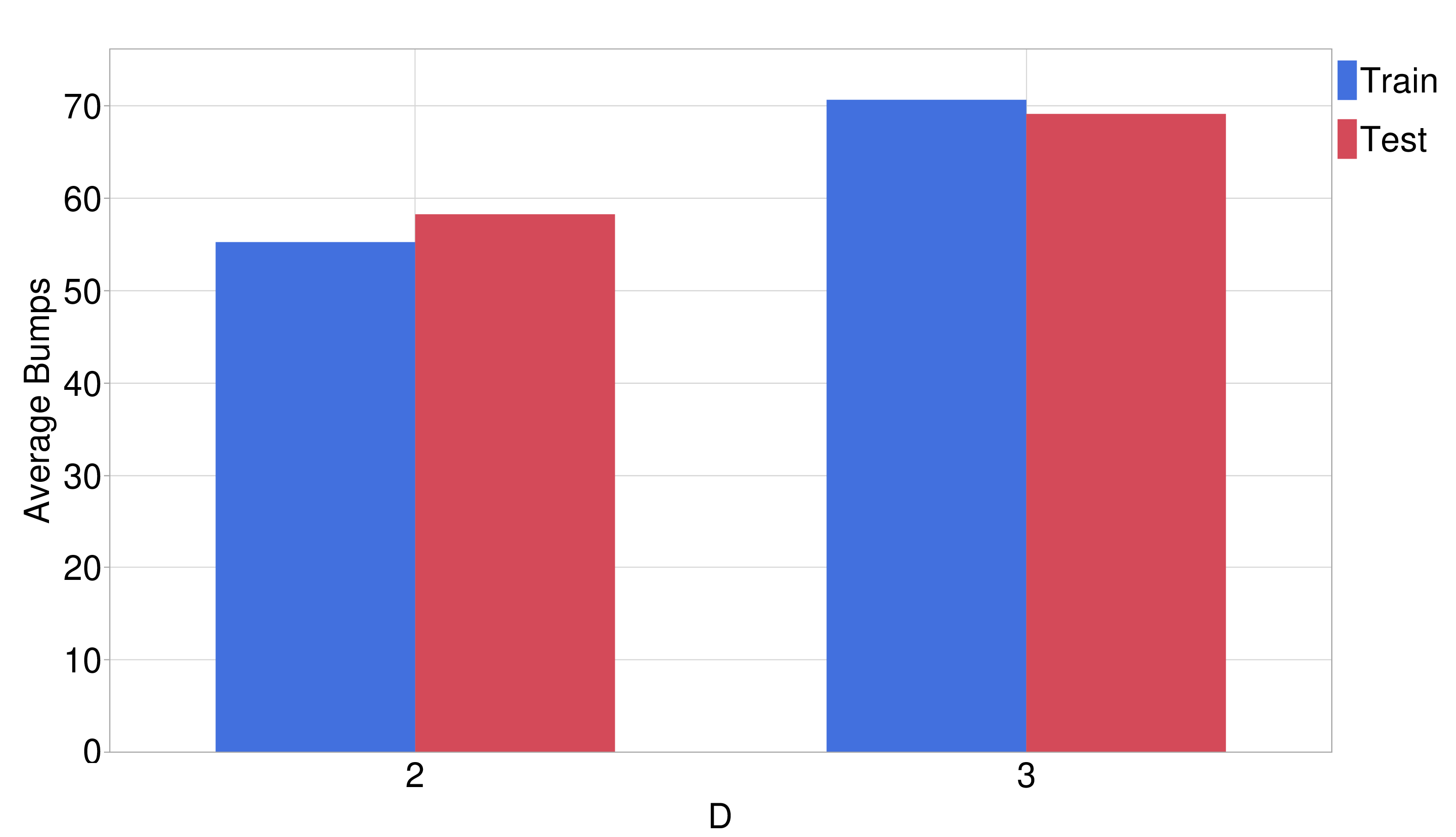}
    \label{stochfigure1} 
     \end{minipage}%
}
\resizebox{7.5cm}{!}{
    \begin{minipage}{.48\textwidth}
        \centering
\caption{Train vs Test performance of $\textbf{\textsc{MIP}}_{\textbf{DNTP-S}}$ - Average Vacant shifts}

\includegraphics[width=\textwidth]{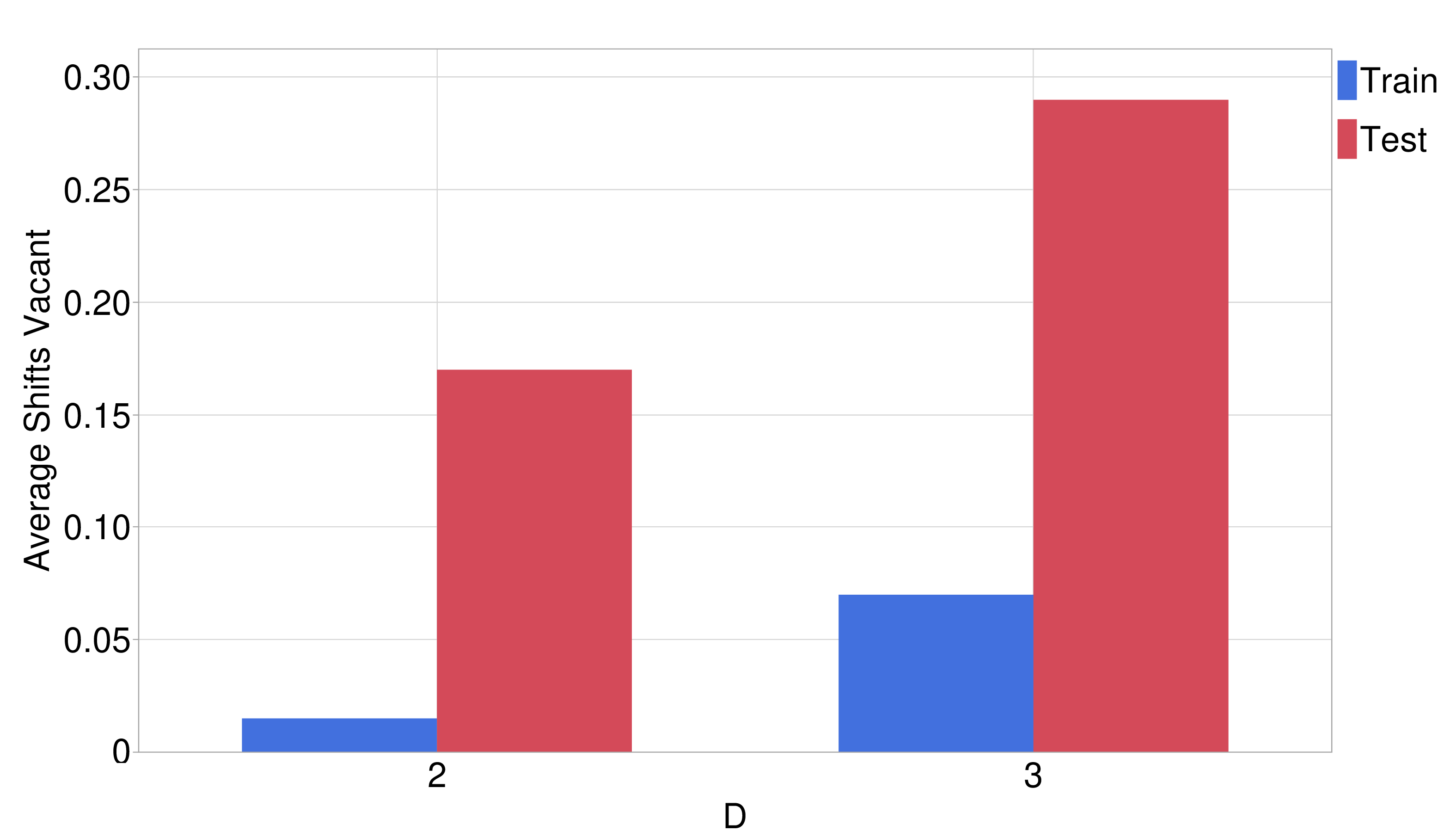}       \label{stochfigure2}
    \end{minipage}
    }
\end{figure}

First, we analyze the performance of a stochastic policy by solving ${\textbf{\textsc{MIP}}_{\textbf{DNTP-S}}}$. We solve the formulation with 200 training scenarios of equal probability for seven days without any acceleration strategy for the solver. The obtained policy, the stochastic policy, is then simulated over 500 test instances. Each vacancy is penalized with a cost of $G = 200$. The results are given in \Cref{stochfigure1} and \Cref{stochfigure2}. The average bumps are similar for both the train and test instances. However, a major difference is seen in the vacancies. The average vacancies increase considerably in the test instances compared to the train instances. They do not meet the 0.3\% threshold indicated by the orange line in the graph. There seems to be an opportunity to improve on the number of vacant shifts.

Additionally, we present results in the E-Companion that compare the simulation outcomes of the policy derived from this formulation with those of the PFAs across 200 train instances. The main takeaway from this comparison is that the ONP policy can closely approximate the solution produced by the stochastic formulation. Readers are encouraged to see these plots.

 \subsubsection{Comparing Stochastic Policy and PFA on test instances}  The decision-maker has two criteria to work with: reducing total bumps and allocating all shifts to employees. To select a policy to use with the electronic system, the management can establish weights for each criterion and choose a policy with minimum cost. Our industrial partner would also like a policy with a maximum shift vacancy of $0.3\%$, which corresponds to an average vacancy of 0.15 shifts. This vacancy percentage is tolerated because some employees may respond after the horizon ends and some manual changes are made, and the management is confident in scheduling this fraction of unscheduled shifts.  Also, it is possible to use overtime. However, relying on these options is generally costly. \Cref{tab:res4} highlights the best parameters found using the 500 validation instances. For offline solution-based policies, high percentile levels work very well. The intuition here is that for a single instance, the offline solution has information about the future and can decide exactly when to notify someone and when to wait. However, in several instances, one needs to protect against vacant shifts as they have a very high penalty. The NAW policy with parameters $\eta = 3$ and $w = 7$ yielded the best results in our tests. In contrast, the company had been using parameters $\eta = 5$ and $w = 1$, which led to a higher number of bumps due to the increased frequency of notifications per minute.

 \Cref{tab:res5} gives the performance of these policies with the best-found parameters in terms of total bumps for 500 test instances. ONP works best (indicated by bold) for both the cutoff of 2 and 3 hours, respectively, in terms of cost as compared to NAW policy. ONP has significantly fewer bumps, slightly compromising on vacancies from  \Cref{tab:res7}. Even the total cost is the least for ONP in \Cref{tab:res6}. Note that meeting the shift vacancy constraint is more important. The policy obtained from $\textbf{\textsc{MIP}}_{\textbf{DNTP-S}}$ is very poor regarding shift vacancies. This policy does not meet the management requirements of 0.3\% average vacancies highlighted by the red color in \Cref{tab:res5}. Thus, this shows that a dynamic threshold policy is better than the current static policy. The NA policy shows a huge number of bumps, thus highlighting the shortcomings of a greedy approach. 
  
\begin{table}
\centering
\resizebox{6cm}{!}{
\parbox[t][][t]{.47\linewidth}{
\centering
\caption{Best Policy Parameters that minimize average bumps subject to 0.3\% shift vacancy } 
\label{tab:res4}
\centering
\resizebox{\linewidth}{!}{
    {\def\arraystretch{1.2}  
\begin{tabular}{|c|cccc|}
\hline
\textbf{}  & \multicolumn{4}{c|}{\textbf{Policy}}                                                                                                                   \\ \hline
\textbf{D} & \multicolumn{1}{c|}{\textbf{NA}}  & \multicolumn{1}{c|}{\textbf{NAW}}  & \multicolumn{1}{c|}{\textbf{ONP}} & $\textbf{\textsc{MIP}}_{\textbf{DNTP-S}}$ \\ \hline
\textbf{2} & \multicolumn{1}{c|}{150} & \multicolumn{1}{c|}{3, 7}   & \multicolumn{1}{c|}{95}    & -           \\ \hline
\textbf{3} & \multicolumn{1}{c|}{150} & \multicolumn{1}{c|}{3, 7}   & \multicolumn{1}{c|}{95}    & -           \\ \hline
\end{tabular}}}

}}
\hfill
\centering
\resizebox{6cm}{!}{
\parbox[t][][t]{.47\linewidth}{
\centering
\caption{Average Bumps for the best policy on test instances } 
\label{tab:res5}
\centering
\resizebox{\linewidth}{!}{
    {\def\arraystretch{1.3}  
\begin{tabular}{|c|cccc|}
\hline
\textbf{}  & \multicolumn{4}{c|}{\textbf{Policy}}                                                                                                                                                                       \\ \hline
\textbf{D} & \multicolumn{1}{c|}{\textbf{NA}} & \multicolumn{1}{c|}{\textbf{NAW}} &  \multicolumn{1}{c|}{\textbf{ONP}} &  $\textbf{\textsc{MIP}}_{\textbf{DNTP-S}}$ \\ \hline
\textbf{2} & \multicolumn{1}{c|}{618.90}      & \multicolumn{1}{c|}{83.13}            & \multicolumn{1}{c|}{\textbf{69.42}}          & \textcolor{red}{58.28}                \\ \hline
\textbf{3} & \multicolumn{1}{c|}{671.3}      & \multicolumn{1}{c|}{115.80}           & \multicolumn{1}{c|}{\textbf{82.63}}          & \textcolor{red}{69.15}             \\ \hline
\end{tabular}}}

}}

\end{table}

\begin{table}
\centering
\resizebox{6cm}{!}{
\parbox[t][][t]{.47\linewidth}{
\caption{Average Cost for the best policy on test instances  } 
\label{tab:res6}

\resizebox{\linewidth}{!}{
    {\def\arraystretch{1.3}  
\begin{tabular}{|c|cccc|}
\hline
\textbf{}  & \multicolumn{4}{c|}{\textbf{Policy}}                                                                                                                   \\ \hline
\textbf{D} & \multicolumn{1}{c|}{\textbf{NA}}  & \multicolumn{1}{c|}{\textbf{NAW}}  & \multicolumn{1}{c|}{\textbf{ONP}} & $\textbf{\textsc{MIP}}_{\textbf{DNTP-S}}$ \\ \hline
\textbf{2} & \multicolumn{1}{c|}{618.90}      & \multicolumn{1}{c|}{88.37}            & \multicolumn{1}{c|}{\textbf{83.42}}          & \textcolor{red}{91.88}\\ \hline
\textbf{3} & \multicolumn{1}{c|}{671.3}      & \multicolumn{1}{c|}{121.00}           & \multicolumn{1}{c|}{\textbf{102.23}}          & \textcolor{red}{127.55}\\ \hline
\end{tabular}}}

}
}
\hfill
\resizebox{6cm}{!}{
\parbox[t][][t]{.47\linewidth}{
\caption{Average Shift Vacancies for the best policy  on test instances} 
\label{tab:res7}
\centering
\resizebox{\linewidth}{!}{
    {\def\arraystretch{1.17}  
\begin{tabular}{|c|cccc|}
\hline
\textbf{}  & \multicolumn{4}{c|}{\textbf{Policy}}                                                                                                                                                                       \\ \hline
\textbf{D} & \multicolumn{1}{c|}{\textbf{NA}} & \multicolumn{1}{c|}{\textbf{NAW}} &  \multicolumn{1}{c|}{\textbf{ONP}} &  $\textbf{\textsc{MIP}}_{\textbf{DNTP-S}}$ \\ \hline
\textbf{2} & \multicolumn{1}{c|}{0}      & \multicolumn{1}{c|}{0.026}            & \multicolumn{1}{c|}{\textbf{0.07}}          & \textcolor{red}{0.17}                \\ \hline
\textbf{3} & \multicolumn{1}{c|}{0}      & \multicolumn{1}{c|}{0.025}           & \multicolumn{1}{c|}{\textbf{0.098}}          & \textcolor{red}{0.29}             \\ \hline
\end{tabular}}}

}}

\end{table}

To ensure fairness in comparing the stochastic formulation, we train both the stochastic policy and the ONP using the same number of training instances (200). The results are summarized in \Cref{tab:stochvsONP}. First, in terms of test and train performance for both cutoff values, the test performance is worse than the train performance for both policies. Second, while there is minimal difference in the number of bumps between the training and testing instances for both policies, we observe a notable increase in vacancy in the test cases. This increase is more pronounced for the stochastic policy, which fails to meet the 0.3\% vacancy criterion. The ONP, on the other hand, appears to generalize better to unseen instances. This indicates that ONP policies perform better, even on small training sets. Additionally, when comparing ONP trained on 1000 and 200 instances, precisely the bumps and vacancies in \Cref{tab:stochvsONP}, \Cref{tab:res5}, \Cref{tab:res6}, and \Cref{tab:res7}, the results with 1000 training instances are slightly better than those trained on 200 instances. 
\begin{table}[!htbp]

\centering
\caption{Comparing Stochastic Policy and ONP on 200 train and 500 test instances }
\label{tab:stochvsONP}
\resizebox{14cm}{!} 
{ 
    {\def\arraystretch{1.25}  
\begin{tabular}{|c|ccccccc|}
\hline
\textbf{}                                     & \multicolumn{7}{c|}{\textbf{Policy}}                                                                                                                                                                                                             \\ \hline
\textbf{}                                     & \multicolumn{1}{c|}{\textbf{}}  & \multicolumn{3}{c|}{${\textbf{\textsc{MIP}}_{\textbf{DNTP-S}}} $}                                                                          & \multicolumn{3}{c|}{\textbf{ONP}}                                                           \\ \hline
\textbf{}                                     & \multicolumn{1}{c|}{\textbf{D}} & \multicolumn{1}{c|}{\textbf{Cost}} & \multicolumn{1}{c|}{\textbf{Bumps}} & \multicolumn{1}{c|}{\textbf{Vacant Shifts}} & \multicolumn{1}{c|}{\textbf{Cost}} & \multicolumn{1}{c|}{\textbf{Bumps}} & \textbf{Vacant Shifts} \\ \hline
\multirow{2}{*}{\textbf{200 Train Instances}} & \multicolumn{1}{c|}{\textbf{2}} & \multicolumn{1}{c|}{\textbf{58.26}}         & \multicolumn{1}{c|}{55.26}          & \multicolumn{1}{c|}{0.015}            & \multicolumn{1}{c|}{73.95}         & \multicolumn{1}{c|}{69.95}          & 0.02            \\ \cline{2-8} 
                                              & \multicolumn{1}{c|}{\textbf{3}} & \multicolumn{1}{c|}{\textbf{84.68}}         & \multicolumn{1}{c|}{70.68}          & \multicolumn{1}{c|}{0.07}             & \multicolumn{1}{c|}{96.15}         & \multicolumn{1}{c|}{87.15}          & 0.045            \\ \hline
\multirow{2}{*}{\textbf{500 Test Instances}}  & \multicolumn{1}{c|}{\textbf{2}} & \multicolumn{1}{c|}{92.28}         & \multicolumn{1}{c|}{58.28}          & \multicolumn{1}{c|}{\textcolor{red}{0.17}}             & \multicolumn{1}{c|}{{86.80}}         & \multicolumn{1}{c|}{\textbf{70.80}}          & 0.08             \\ \cline{2-8} 
                                              & \multicolumn{1}{c|}{\textbf{3}} & \multicolumn{1}{c|}{127.55}        & \multicolumn{1}{c|}{69.15}          & \multicolumn{1}{c|}{\textcolor{red}{0.29}}             & \multicolumn{1}{c|}{{105.11}}        & \multicolumn{1}{c|}{\textbf{82.71}}          & 0.11             \\ \hline
\end{tabular}} }

\end{table}

\subsection{Policy Evaluation - Different Preferences} 

One expects that each employee will have different shift preferences in real life. The preferences can be cognitive, cultural, social, temporal, situational, health conscious, or based on familiarity (\cite{rogers2003diffusion}, \cite{daniel2017thinking}). Without real-world preferences, we test some synthetic cases and see how the offline solution-based policies compare with the other policies. Let $\mathcal{L} = \{1,\dots,L\}$ denote the set of shifts. Also,  let $\mathcal{X}_i = \{p_1,\dots,p_L\}$ denote the preference set of employee $i$. Thus, if $l = p_o$, it denotes that shift $l$ occurs in position $o \in \{1,\dots, L\}$ of the preference set. We consider the following preferences:
\begin{itemize}
\item \textbf{Fixed Ranked Preferences} - This distribution represents the base case where all preference sets are identical.
\item \textbf{Undesirable Preferences} -  Certain employees may dislike certain shifts. We assume that each employee randomly disapproves of $\Bar{l}$ shifts from the full set. Starting from a fixed ranked preference set, we randomly move $\Bar{l}$ shifts to the end of the preference list. These undesirable shifts remain consistent for any employee across all simulations.
    \item \textbf{Grouped Preferences} - Shifts are assumed to be divided into two groups of equal sizes. Any employee strictly prefers shifts from one set over the other. Within a group, all shifts follow the same preference order for all employees. If the shifts from the preferred group are exhausted, the employee can choose from the other group. The preference distribution represents where shifts can be divided into day and night shifts, and certain employees prefer one over the other. Each employee chooses a group with the same probability.  
    \item \textbf{Perturbed Preferences} - We assume that there is still a general fixed ranked order of the shifts from 0 to $L$. However, the shift can deviate slightly from the order for a given employee. Specifically, any shift $l $ can occur in position $o \in \{l - 4,\dots,l + 4\}$ with equal probability. In other words, we have a sliding window around the original shift position for the shift.
    \item \textbf{Perturbed Preferences with Undesirable shifts} - This represents a combination of perturbed preferences and undesirable shifts. 
    \item \textbf{Uniform Preferences} - Employees are equally likely to choose any shift from those available. This represents a very unrealistic scenario, but we use it as a baseline scenario.
\end{itemize}

The best parameters for each preference distribution across all policies are provided in \Cref{tab:res8}, subject to the maximum vacancy constraint of 0.3\%. The various preferences are listed roughly in the order of randomness observed, starting from fixed-ranked preferences. Consequently, we observe that more employees can be notified across all policies, as the best parameters generally increase for the same policy. \Cref{tab:res9} presents all policies' average number of bumps, subject to the shift vacancy constraint. Note that all bumps, in this case, represent realized bumps. The ONP continued to perform the best across all preference distributions. However, as randomness increases, the performance gain over NAW diminishes. If employees behaved randomly (Uniform), the average number of bumps was very low. In such a scenario, one would notify the maximum number of employees, contrasting with the worst-case scenario of the same preferences. The table suggests that our designed policy will result in fewer bumps than the current NAW policy. Overall, ONP still proves to be very robust even with different preferences, suggesting such a policy may be fruitful if implemented.

\begin{table}[ht!]
\centering
\caption{Best Policy Parameters that minimize bumps in average subject to 0.3\% shift vacancy on validation instances for real-world data} 
\label{tab:res8}
\resizebox{9cm}{!} 
{ 
    {\def\arraystretch{0.9}  
    
\begin{tabular}{|ccccc|}
\hline
\multicolumn{1}{|c|}{\textbf{}}                   & \multicolumn{1}{c|}{\textbf{}}                    & \multicolumn{3}{c|}{\textbf{Policy}}                                                                                                                                     \\ \hline
\multicolumn{1}{|c|}{\textbf{D}}                  & \multicolumn{1}{c|}{\textbf{Preference}}          & \multicolumn{1}{c|}{\textbf{NA}} & \multicolumn{1}{c|}{\textbf{NAW}} & \multicolumn{1}{c|}{\textbf{ONP}} \\ \hline
\multicolumn{1}{|c|}{\multirow{4}{*}{\textbf{2}}} & \multicolumn{1}{c|}{\textbf{Fixed Ranked}}        & \multicolumn{1}{c|}{150}         & \multicolumn{1}{c|}{3, 7}         & \multicolumn{1}{c|}{95}           \\ \cline{2-5} 
\multicolumn{1}{|c|}{}                            & \multicolumn{1}{c|}{\textbf{Undesirable}}         & \multicolumn{1}{c|}{150}         & \multicolumn{1}{c|}{3, 7}         & \multicolumn{1}{c|}{95}           \\ \cline{2-5} 
\multicolumn{1}{|c|}{}                            & \multicolumn{1}{c|}{\textbf{Grouped}}             & \multicolumn{1}{c|}{150}         & \multicolumn{1}{c|}{2, 4}         & \multicolumn{1}{c|}{95}           \\ \cline{2-5} 
\multicolumn{1}{|c|}{}                            & \multicolumn{1}{c|}{\textbf{Perturbed}}           & \multicolumn{1}{c|}{150}         & \multicolumn{1}{c|}{2, 4}         & \multicolumn{1}{c|}{98}           \\ \cline{2-5} 
\multicolumn{1}{|c|}{}                            & \multicolumn{1}{c|}{\textbf{Per. w. Undesirable}} & \multicolumn{1}{c|}{150}         & \multicolumn{1}{c|}{2, 4}         & \multicolumn{1}{c|}{98}           \\ \cline{2-5} 
\multicolumn{1}{|c|}{}                            & \multicolumn{1}{c|}{\textbf{Uniform}}             & \multicolumn{1}{c|}{150}         & \multicolumn{1}{c|}{3, 6}         & \multicolumn{1}{c|}{98}           \\ \hline
\multicolumn{5}{|c|}{}                                                                                                                                                                                            \\ \hline
\multicolumn{1}{|c|}{\multirow{4}{*}{\textbf{3}}} & \multicolumn{1}{c|}{\textbf{Fixed Ranked}}        & \multicolumn{1}{c|}{150}         & \multicolumn{1}{c|}{3, 7}         & \multicolumn{1}{c|}{95}           \\ \cline{2-5} 
\multicolumn{1}{|c|}{}                            & \multicolumn{1}{c|}{\textbf{Undesirable}}         & \multicolumn{1}{c|}{150}         & \multicolumn{1}{c|}{3, 7}         & \multicolumn{1}{c|}{98}           \\ \cline{2-5} 
\multicolumn{1}{|c|}{}                            & \multicolumn{1}{c|}{\textbf{Grouped}}             & \multicolumn{1}{c|}{150}         & \multicolumn{1}{c|}{1, 2}         & \multicolumn{1}{c|}{98}           \\ \cline{2-5} 
\multicolumn{1}{|c|}{}                            & \multicolumn{1}{c|}{\textbf{Perturbed}}           & \multicolumn{1}{c|}{150}         & \multicolumn{1}{c|}{2, 4}         & \multicolumn{1}{c|}{98}           \\ \cline{2-5} 
\multicolumn{1}{|c|}{}                            & \multicolumn{1}{c|}{\textbf{Per. w. Undesirable}} & \multicolumn{1}{c|}{150}         & \multicolumn{1}{c|}{4, 8}         & \multicolumn{1}{c|}{98}           \\ \cline{2-5} 
\multicolumn{1}{|c|}{}                            & \multicolumn{1}{c|}{\textbf{Uniform}}             & \multicolumn{1}{c|}{150}         & \multicolumn{1}{c|}{4, 8}         & \multicolumn{1}{c|}{98}           \\ \hline
\end{tabular}}}

\end{table}
\begin{table}
\centering
\caption{Average bumps for each policy on real data-based test instances across various preference distributions}
\label{tab:res9}
\resizebox{9cm}{!} 
{ 
    {\def\arraystretch{1.25}  
\begin{tabular}{|cccccc|}
\hline
\multicolumn{1}{|c|}{\textbf{}}                   & \multicolumn{1}{c|}{\textbf{}}                    & \multicolumn{4}{c|}{\textbf{Policy}}                                                                                                                                                                       \\ \hline
\multicolumn{1}{|c|}{\textbf{D}}                  & \multicolumn{1}{c|}{\textbf{Preference}}          & \multicolumn{1}{c|}{\textbf{NA}} & \multicolumn{1}{c|}{\textbf{NAW}}  & \multicolumn{1}{c|}{\textbf{ONP}} &  \textbf{$\textbf{\textsc{MIP}}_{\textbf{DNTP-S}}$} \\ \hline
\multicolumn{1}{|c|}{\multirow{4}{*}{\textbf{2}}} & \multicolumn{1}{c|}{\textbf{Fixed Ranked}}               & \multicolumn{1}{c|}{623.73}      & \multicolumn{1}{c|}{83.86}            & \multicolumn{1}{c|}{\textbf{70.85}}               & \textcolor{red}{54.56}               \\ \cline{2-6} 
\multicolumn{1}{|c|}{}                            & \multicolumn{1}{c|}{\textbf{Undesirable}}         & \multicolumn{1}{c|}{424.65}      & \multicolumn{1}{c|}{50.20}             & \multicolumn{1}{c|}{\textbf{44.50}}               & \textcolor{red}{37.42}                \\ \cline{2-6} 
\multicolumn{1}{|c|}{}       & \multicolumn{1}{c|}{\textbf{Grouped}}             & \multicolumn{1}{c|}{347.05}      & \multicolumn{1}{c|}{43.25}       
& \multicolumn{1}{c|}{\textbf{33.93}}              & \textcolor{red}{29.28}            \\ \cline{2-6}                   
\multicolumn{1}{|c|}{}                            & \multicolumn{1}{c|}{\textbf{Perturbed}}             & \multicolumn{1}{c|}{225.46}      & \multicolumn{1}{c|}{28.93}          & \multicolumn{1}{c|}{\textbf{27.35}}             & \textcolor{red}{19.89}            \\ \cline{2-6} 
\multicolumn{1}{|c|}{}                            & \multicolumn{1}{c|}{\textbf{Per. w. Undesirable}} & \multicolumn{1}{c|}{197.50}      & \multicolumn{1}{c|}{25.30}               & \multicolumn{1}{c|}{\textbf{23.78} }              & \textcolor{red}{17.21}            \\ \cline{2-6} 
\multicolumn{1}{|c|}{}                            & \multicolumn{1}{c|}{\textbf{Uniform}}             & \multicolumn{1}{c|}{57.71}       & \multicolumn{1}{c|}{5.23}             & \multicolumn{1}{c|}{\textbf{4.12}  }             & \textcolor{red}{4.20}            \\ \hline
\multicolumn{6}{|l|}{}                                                                                                                                                                                                                                                                                             \\ \hline
\multicolumn{1}{|c|}{\multirow{4}{*}{\textbf{3}}} & \multicolumn{1}{c|}{\textbf{Fixed Ranked}}               & \multicolumn{1}{c|}{678.97}      & \multicolumn{1}{c|}{118.32}              & \multicolumn{1}{c|}{\textbf{84.84}}           & \textcolor{red}{64.87}               \\ \cline{2-6} 
\multicolumn{1}{|c|}{}                            & \multicolumn{1}{c|}{\textbf{Undesirable}}         & \multicolumn{1}{c|}{464.5}       & \multicolumn{1}{c|}{67.72}               & \multicolumn{1}{c|}{\textbf{63.00}   }            & \textcolor{red}{44.27}              \\ \cline{2-6} 
\multicolumn{1}{|c|}{}             & \multicolumn{1}{c|}{\textbf{Grouped}}             & \multicolumn{1}{c|}{380.47}      & \multicolumn{1}{c|}{56.08}            & \multicolumn{1}{c|}{\textbf{48.24}}               & \textcolor{red}{35.09}              \\ \cline{2-6}                
\multicolumn{1}{|c|}{}                            & \multicolumn{1}{c|}{\textbf{Perturbed}}             & \multicolumn{1}{c|}{246.73}      & \multicolumn{1}{c|}{37.74}         & \multicolumn{1}{c|}{\textbf{32.69}}           & \textcolor{red}{23.16}                \\ \cline{2-6} 
\multicolumn{1}{|c|}{}                            & \multicolumn{1}{c|}{\textbf{Per. w. Undesirable}} & \multicolumn{1}{c|}{216.41}      & \multicolumn{1}{c|}{33.44}             & \multicolumn{1}{c|}{\textbf{28.25} }                & \textcolor{red}{20.35}                \\ \cline{2-6} 
\multicolumn{1}{|c|}{}                            & \multicolumn{1}{c|}{\textbf{Uniform}}             & \multicolumn{1}{c|}{65.30}       & \multicolumn{1}{c|}{6.45}              & \multicolumn{1}{c|}{\textbf{5.18}  }               & \textcolor{red}{5.17}                 \\ \hline
\end{tabular}}}

\end{table}

\subsection{Managerial Insights}
The analysis has provided several important insights for decision-makers. First, since preferences can be difficult to model, decision-makers should aim to minimize potential bumps, which will, in turn, reduce actual bumps. Second, the results from offline solutions demonstrate that the offline Notification Timing Problem (\textsc{NTP}) is challenging; however, with complete information, potential bumps can be nearly eliminated. Third, the current policy lacks the flexibility to adapt to the system's evolving states. To address this issue, we developed a dynamic threshold policy (ONP) based on the states observed in the offline solution. Despite its simple structure, our results show that the ONP significantly reduces both bumps and shift vacancies compared to the most tuned version of the current policy (NAW). This underscores the effectiveness of more advanced policies and suggests that managers should adopt such policies to improve scheduling outcomes. Furthermore, through \Cref{alg:two}, we provide management with an easy-to-implement methodology for defining and tuning various threshold-based policies. This allows management to experiment with different tuned policies and select the most effective one.

Analysis of real-world data revealed that only 20\% of employees respond promptly (within 5 minutes) to shift notifications, while about 50\% do not respond at all. If such early and slow responders could be identified, they could be targeted by adapting the ONP. Systems that can effectively segment employees by response time and tailor the notification strategy will see improved shift-filling rates. ONP remains scalable and can be redone without considerable effort as it is compiled offline.


%% file: conclusion.tex
\section{Conclusion}\label{conclusion}

Traditional on-demand personnel scheduling systems face significant challenges regarding employee flexibility and operational ease for management. Despite their widespread use, these systems have received limited attention in the operations research community for optimization purposes. This paper introduces a novel, flexible, dynamic personnel scheduling system to tackle these challenges. A key feature of our system is the ability for senior employees to replace or bump junior employees already assigned to shifts. These bumps are undesirable, and management aims to minimize them while ensuring all shifts are filled. Given the difficulty of capturing employee preferences, we consider the worst-case scenario of minimizing potential bumps.
Our main contribution is demonstrating that the problem of reducing bumps while filling all shifts is $\mathcal{NP}$-complete, even when all uncertainty is revealed in advance, through a reduction from the \textsc{SUBSET-SUM} problem. Furthermore, we show that any online algorithm can perform arbitrarily poorly in certain instances of the dynamic version of the problem. 
To address these issues, we propose a Policy Function Approximation (PFA) with a threshold structure for the dynamic version. Our approach estimates the threshold by heuristically solving a 2 stage formulation of the dynamic problem.  The methodology solves several individual instances offline and then aggregates the results from all the instances using simple descriptive statistics at each decision epoch. These threshold estimates are validated on unseen cases. The policy is developed under the worst-case assumption of uniform preferences, demonstrating superior performance compared to existing methods when applied to both uniform and synthetic preference distributions.

Several future extensions for the problem are possible. First, we presented a simpler system version that only operates over one day. In reality, the scheduling can start four days before a list of shifts starts. Second, the cutoff time to bump someone could be a function of the time remaining and decrease as the system moves forward into the planning period. Next, the number of shifts available to schedule varies, and the management would like a policy to consider this. Our industrial partner also has additional employee data that could be used to fine-tune the policies. We plan to tackle these challenges in future research. 

%% file: appendix.tex
%

\section{Pseudo Code to find a bump chain}

Let $I = \{I_l\}_{l \in \mathcal{L}}$ and $L = \{L_j\}_{j \in \mathcal{E}}$. 

\begin{algorithm}[!ht]
   \SetKwInOut{Input}{Input}
    \SetKwInOut{Output}{Output}
    \SetKwFunction{Function}{Compile}
    \SetKwProg{Fn}{Function}{: }{end}
\caption{To construct $\mathcal{B}_i$}
\label{alg:three}
\Fn{\Function{$i, e_i, I, L$}}{  

$\mathcal{B}_i = \{\}$\;
$i^\prime \gets i$\;
\While{$I_{L_{i^\prime}} \neq \phi$}{

    $\mathcal{B}_i = \mathcal{B}_i \cup \{I_{L_{i^\prime}}\}$\;
 $i^\prime \gets I_{L_{i^\prime}}$\; 
 }

  $i \gets i + 1$\;
  
    \Return $\mathcal{B}_i$;
  }
\end{algorithm}

\section{Proofs}
\label{append}
Note: Any reference to a bump is a reference to a potential bump.
\begin{repeattheorem}
Given an instance $I = \{H, M, \mathbf{r}\}$, let $S^*_\mathcal{X}$  denote an optimal schedule with preferences $\mathcal{X}$ such that $\mathcal{X} = \{{\mathcal{X}_i}\}_{i \in \mathcal{E}}$. Also, let $S^*_p$ be the schedule that minimizes the total potential bumps for this instance, and  $S^*_I$ be the schedule that minimizes the total bumps with identical preferences. Then, $B_{S^*_\mathcal{X}} \leq B_{S^*_p} = B_{S^*_I} $.
\end{repeattheorem}
\proof{\textbf{Proof of Theorem 1}} We show this in two parts. 

\begin{itemize}
    \item[1)]  It is trivial to say that $B_{S^*_\mathcal{X}} \leq B_{S^*_p}$. This is because potential bumps are all cases when a senior employee has responded later than a junior employee. 
    \item[2)]  Now consider the case when preferences are identical. By 1), $B_{S^*_I} \leq B_{S^*_p}$. Let $\{1, \dots, M\}$ represent the set of shifts where 1 represents the most preferred shift, 2 represents the second most preferred and so on in the identical preferences setting. Consider an employee $i$ responding back at $e_i$ and choosing the current preferred and available shift $L_i$. Let 
    \begin{align}
         \mathcal{P}_{i} = \{j_1, j_2, \dots, j_k : e_{j_o} < e_i, i<{j_1}<\dots<j_k\}
    \end{align}
    
    represent the potential bumps as in Equation 2.  $j_1$ represents the least junior employee responding before $i$; ($j_1 = \{\min j: j>i, e_i > e_j\, j \in \mathcal{E}\}$) and is the first potential candidate to get bumped if $i$ responds. If $\mathcal{B}_{i} = \phi$, there are no potential bumps induced by a response from $i$.
    
Consider any $i$ such that $\mathcal{P}_{i} \neq \phi$, that is,  there are some potential bumps.  Under identical preferences, $\mathcal{X}_i = \mathcal{X}_j = \{1,\dots, M\}, \forall j \in \mathcal{P}_i$. We now show that $\mathcal{P}_{i}$ is also a bump chain defined by Equation 1, for that employee $i$ when preferences are identical. That is, all potential bumps are actually realised bumps when preferences are identical ($\mathcal{P}_i = \mathcal{B}_i$) for all employees $i \in \mathcal{E}$ and, as a consequence, $B_{S^*_I} \geq B_{S^*_p}$.
  
    Let $I_{l} = j_1$ at time $e_i$, that is employee $j_1$ occupies shift $l$.  All shifts in the set $\{1,\dots,l-1\}$ are chosen by senior employees to $j_1$, otherwise $j_1$ would occupy them as they more preferred in $\mathcal{X}_{j_{1}}$. This also means that $L_j = l + 1$. Consequently, for any employee in $j_o \in\mathcal{P}_i$, $I_{l+o-1} = j_{o}$. Also

    \begin{align}
L_{j_o} = 
\begin{cases}
l + o & \text{if } l + o \leq M ;\\
\phi & \text{else;}\\
\end{cases}
\label{nextshift}
\end{align}

    As $i$ is senior to $j_1$, shift $l$ will always be available to him as long as $j_1$ occupies it.  When $i$ responds,  he too will not have any shifts available to him in the set $\{1,\dots,l-1\}$. If such a shift were available, then $j$ would occupy it because it is more preferred. Hence $l = L_i$ and $I_l = j_1$ at the time $i$ responds. Employee $i$ will bump $j_1$, and $j_1$ is the first element of the bump set in Equation 1.  As a consequence, $j_1$ will then similarly bump $j_2$ for the next shift $l+1$ as $L_j = l + 1$. Thus, this will induce a chain of bumps for all employees in $\mathcal{P}_i$. For the last person in the chain, either there is no shift available for $j_1$, or his preferred shift will be unoccupied from \Cref{nextshift}. Thus, all employees in $\mathcal{P}_{i}$ are bumped $\forall i \in \mathcal{E}$, and it is a bump chain $\mathcal{B}_i$.  Therefore, $B_{S^*_p} \leq B_{S^*_I}$ and $B_{S^*_p} = B_{S^*_I} $. 
\end{itemize}

\begin{repeatproposition}
    The minimum makespan of a NBS is $C^*_0 = r_1 + \sum\limits_{i=1}^{n} {(r_{i+1} - r_i)^+}$.
\end{repeatproposition}

\proof{\textbf{Proof of Proposition 1}} 
\textbf{Case~1}: $r_i \leq r_{i+1}, \forall i \in \mathcal{E}$. Therefore, we can send notifications to employees $i$ and $i + 1$ at the same time ($s_i = s_{i+1}$) without introducing a bump. The makespan of the schedule is increased by  $r_{i+1} - r_i $. \\
\textbf{Case~2}: $r_i > r_{i+1}, \forall i \in \mathcal{E}$.
The only way to avoid a bump, in this case, is to make sure that both of them respond at the same time $e_i = e_{i+1}$. This is easily ensured by setting $s_{i+1} = s_i + r_i - r_{i+1}$. Hence, the makespan is not increased in this case.

\endproof

\begin{repeatproposition}
 The block schedule has the following properties:
\begin{enumerate}
    \item $\forall k \in A, s_{i_k} = s_{i_k - 1}$.
    \item $\forall i \in \mathcal{E}^S_{k}, \forall k \in A$,  $e_i = e_{i_k}$ or $s_i = s_{i_k}$.
    \item $s_{i_{k+1}} = s_{i_{k}} + \mathbf{1}_{k \notin A^\prime}a_{k}$, where $\mathbf{1}$ denotes indicator variable.
    \item $B(S(A^\prime)) = \sum\limits_{j \in A^\prime} a_j$
    
    \item $C(S(A^\prime)) = C^*_0 - \sum\limits_{j \in A^\prime } a_j $.
\end{enumerate}
\end{repeatproposition}
\proof{\textbf{Proof of Proposition 2}}
Property 1, 2, and 3 are true by construction.  
\begin{itemize}
    \item[4.] Consider some $k \in A^\prime$, $r_{i_k} > r_i$ from Equation (7) in the main article, and $s_{i_k} = s_i$ from Equation (8) in the main article, $\forall i \in \mathcal{E}^S_{k}$. As a consequence  
  \begin{equation*}
        e_{i_k} = r_{i_k} + s_{i_k}  > r_i + s_i = e_i, \forall i_k < i < i_{k+1}. 
  \end{equation*}
Since $a_k = i_{k+1} - i_k$ and using proposition Corollary 1, $b_{i_k} = a_k$.
  \\Now consider some $k \in A\backslash A^\prime$,  and $e_{i_k} = s_{i_k} + r_{i_k}  = s_i + r_i = e_i, \forall  i_k < i < i_{k+1}$. As a consequence, no employee is potentially bumped by $i_k$ and   $b_{i_k} = 0$.
  Hence, the total potential bumps are given by 
  \begin{equation}
  \label{totbumps}
        B(S(A^\prime)) = \sum\limits_{i\in \mathcal{E}} b_{i} = \sum\limits_{k\in A} b_{i_k} = \sum\limits_{k\in A^\prime} a_k 
  \end{equation}

    \item[5.] $C(S(A^\prime)) = C^*_0 - \sum\limits_{j \in A^\prime } a_j $.\\
    Using property 3 and by recursion: 
            \begin{align}
        \label{start_crit}
        s_{i_{k+1}} &= \sum\limits_{j = 1: j\notin A^\prime}^k a_j
    \end{align}
    Now consider the employee with priority $M$. 
\begin{equation}
\label{s_m}
s_{M} =  s_{M-1} =  s_{i_N} + \mathbf{1}_{N \notin A^\prime}a_N.
\end{equation}
Substitute equation \Cref{start_crit} in \Cref{s_m}, 
\begin{align*}
    s_{M} &= \sum\limits_{j = 1: j\notin A^\prime}^{N-1} a_j + \mathbf{1}_{N \notin A^\prime}a_N \nonumber\\
    s_{M} &= \sum\limits_{ j\notin A^\prime} a_j
\end{align*}
As a consequence one can now find $C(S(A^\prime)) = e_M$
\begin{align*}
    C(S(A^\prime)) &= s_M + r_M\\
 &= \sum\limits_{ j\notin A^\prime} a_j + \sum\limits_{ j\in A} a_j\\
     &= 2* \sum\limits_{ j\in A} a_j - \sum\limits_{ j\in A^\prime} a_j \\
     &= C^*_0 - \sum\limits_{ j\in A^\prime} a_j 
\end{align*}

\end{itemize}
    \endproof
\begin{repeatproposition}
For any $\ntpss$, $C = e_M$.
\end{repeatproposition}
\proof{\textbf{Proof of Proposition 3}} 
It suffices to show that $r_M \geq r_i$, since $s_M \geq s_i$, due to seniority $\forall i \in \mathcal{E}$.\\
\textbf{Case 1}: $i = i_k \in \mathcal{E^C}$.
\begin{align*}
    r_{i_k} = \sum\limits_{j = 1}^{k}a_j \leq \sum\limits_{j \in A}a_j = r_M
\end{align*}\\
\textbf{Case 2}:  $ i \in \mathcal{E}^S_{k} : k \in A$.
\begin{align*}
    r_i = r_{i_{k-1}} = \sum\limits_{j = 1}^{k-1}a_j < \sum\limits_{j \in A}a_j = r_M
\end{align*}

Hence $r_M \geq r_i$, $\forall i \in \mathcal{E}$.
\endproof 

\begin{repeatproposition}
$\forall k \in A$, $\forall l \in A$ employee $i_k$ cannot bump employee $i_l$.

\end{repeatproposition}

\proof{\textbf{Proof of Proposition 4}} 
\textbf{\\Case 1} $k \geq l$:  The proposition is true due to seniority $\forall i_k, i_l \in \mathcal{E}$.\\
\textbf{Case 2}  $k < l$ : To show that  $e_{i_k} \leq e_{i_l}.$  
Consider $i_k \in \mathcal{E^C}$ such that $ r_{i_k}  = \sum\limits_{j = 1}^{k}a_j$ from Equation 7 in the main article;
\begin{align*}
    r_{i_{l}} -  r_{i_k} &= \sum\limits_{j = 1}^{l}a_j - \sum\limits_{j = 1}^{k}a_j\\
    r_{i_{l}} -  r_{i_k} &= \sum\limits_{j = k +1 }^{l}a_{j} > 0\\
\end{align*}
Hence $r_{i_{l}} >  r_{i_k}$. Due to seniority $s_{i_{l}} \geq  s_{i_k}$. Adding both inequalities 
\begin{align*}
    r_{i_{l}} + s_{i_{l}} &>  r_{i_k} + s_{i_k}\\
    e_{i_l} &> e_{i_k}
\end{align*}

\endproof
\begin{repeatproposition}
 $\forall i \in \mathcal{E^S}$, $\forall i^\prime \in \mathcal{E}$ employee $i$ cannot bump  employee $i^\prime$.
\end{repeatproposition}

\proof{\textbf{Proof of Proposition 5}} 
\textbf{\\Case 1} $i \geq i^\prime$:  The proposition is true due to seniority $\forall i \in \mathcal{E}^s, i^\prime \in \mathcal{E}$.\\
\textbf{Case 2}  $i < i^\prime$ : To show that  $e_{i} \leq e_{i^\prime}.$  

Response delay for $i$, $ r_{i}  =  r_{i_{k-1}} =  \sum\limits_{j = 1}^{k-1}a_j$ where $i_k< i < i_{k+1}$ and $ k \in A$ from Equation (7) in the main article;\\
\begin{equation}
r_{i^\prime} =   r_{i_{k^\prime - 1 }} : i_{k^\prime} < i^\prime < i_{k^\prime+1}
\end{equation}
Since $i < i^\prime$,  $k \leq k^\prime$ . 

 \begin{align*}
    r_{i^\prime} -  r_{i} &= r_{i_{k^\prime-1 }} -  r_{i_{k-1}} \\
    r_{i^\prime} -  r_{i} &= \sum\limits_{j = 1}^{k^\prime -1}a_j -\sum\limits_{j = 1}^{k-1}a_j\\
    r_{i^\prime} -  r_{i} & \geq 0
\end{align*}
Hence $r_{i^\prime} \geq r_{i}$. Due to seniority $s_{i^\prime} \geq  s_{i}$. As a consequence,  $e_{i^\prime} \geq e_{i}$.
\endproof

\begin{repeatproposition}
$\forall k \in A $, $\forall i  \in \mathcal{E}: i \geq i_{k+1}$, employee $i_k$ cannot bump employee $i$.
\end{repeatproposition}

\proof{\textbf{Proof of Proposition 6}} 
$\forall k \in A $, $\forall i  \in \mathcal{E}: i \geq i_{k+1}$, it suffices to show  $r_{i_k} \leq r_{i}.$  \\
\textbf{Case 1}: For $ i \in \mathcal{E^C} \cap \{i^\prime : i^\prime\geq i_{k+1}, i^\prime \in \mathcal{E}  \} $, the inequality holds due to proposition Proposition 3.\\
\textbf{Case 2}: For $ i \in \mathcal{E^S} \cap \{i^\prime : i^\prime\geq i_{k+1}, i^\prime \in \mathcal{E} \} $, $r_{i} = r_{i_{l-1 }}$ where  $i \in \mathcal{E}^S_k$ and $ l \in A $. Since $i > i_{k+1}$ and $i_{k+1} \in \mathcal{E^C}$,  $l \geq k + 1$.
 \begin{align*}
    r_{i} -  r_{i_k} &= r_{i_{l - 1 }} -  r_{i_k} \\
\end{align*}
Using Proposition 4
, where the response times of stable employees are monotonically non-decreasing, to get $r_{i_{l - 1 }} \geq r_{i_{k + 1 - 1 }} = r_{i_{k  }}$,
 \begin{align*}
     r_{i_{l - 1 }} -  r_{i_k} &\geq r_{i_{k + 1 - 1 }} -  r_{i_k}\\
     r_{i_{l - 1 }} -  r_{i_k} & \geq 0
\end{align*}
Hence $ r_{i} \geq  r_{i_k}$, and since $s_i \geq s_{i_k}$ due to seniority,  $e_i \geq e_{i_k}$ .
\endproof

\begin{repeatcorollary}
$\forall k \in A$, $\forall i\in \mathcal{E}$ if employee $i_k$  bumps employee $i$, then $i \in \mathcal{E}^S_{k}$.
\end{repeatcorollary}

\proof{\textbf{Proof of Corollary 1}} 
This is straightforward from Proposition 3  and Proposition 5   and by construction in Equation (7) from the main article, where $r_{i_k} > r_i $, such that $k \in A$, $i \in \mathcal{E}^S_{k}$. 
\endproof

\begin{repeatlemma}
 $\forall S^* \in \optS\big([a_j]_{j \in \mathcal{A}}, W\big) $, $\forall k \in A, b_{i_k} = a_k $ or $0$. 
\end{repeatlemma}

\proof{\textbf{Proof of Lemma 1}} 
Let $S^* = [s_i, e_i]_{i \in \mathcal{E}} \in \optS\big([a_j]_{j \in \mathcal{A}}, W\big)$, and let $B^* = B(S^*)$ be the optimal number of potential bumps and $C(S^*)$ represent the makespan of the notification schedule. Note we use $B^*, B(S^*)$ interchangeably.
Consider the employee block $\mathcal{E}^S_{k}$, along with critical employee $i_k$, for $k \in A$.  The  following two cases are possible:\\
\textbf{Case 1}: $\forall i $ such that $i \in \mathcal{E}^S_{k}$, $e_i < e_{i_k}$. This directly implies $b_{i_k} = a_k$.
\\
\textbf{Case 2}: $\exists i^\prime$ such that $e_{i^\prime} \geq e_{i_k}$. A new of the schedule $S^\prime = [s^\prime_i, e^\prime_i]_{i \in \mathcal{E}}$ is can be formed by redefining end times for employees of the $\mathcal{E}^S_k$ as follows
\begin{equation*}
e^\prime_i =  \begin{cases}
e_{i_k}& \forall i \in \{i: e_i < e_{i_k}, i \in \mathcal{E}^S_{k}\},  \\
e_i & \text{else}
\end{cases}
\end{equation*}
This new schedule $S^\prime$ is a feasible schedule since $\forall i \in \{i: e_i < e_{i_k}, i \in \mathcal{E}^S_{k}\}, r_i = r_{i_{k - 1}}$ and $e^\prime_{i} =  e_{i_k}$ implies $s^\prime_{i} = s_{i_k} + a_k$. Also $\forall i \in \{i: e_i \geq e_{i_k}, i \in \mathcal{E}^S_{k}\}$,
\begin{align*}
    s^\prime_{i} &= s_i\\
    &  \geq s_{i_{k}} + r_{i_k} - r_{i_{k - 1}}\\
    &  \geq s_{i_{k}} + a_k
\end{align*}
In general, $\forall i \in \mathcal{E}^S_{k}$, $s^\prime_{i} \geq \max\{s_{i_{k}} + a_k , s_{i - 1}\}$.
Hence it is easy to see that for the schedule $S^\prime$, the potential bumps for $i_k$,   $b^\prime_{i_k} = 0$.
 Hence $\forall k^\prime \in  A: k^\prime \neq k$,
\begin{equation}
\label{trans_opt_sched}
    b^\prime_{i_k^\prime} = b_{i_k}
\end{equation} This equation follows from Proposition 6 and Corollary 1. Also from Proposition 5, $b_i = b^\prime_i = 0, \forall i \in \mathcal{E}^S$. Hence for schedule $S^\prime$, 
$b^\prime_{i_k} = 0$ since $\forall i   \in \mathcal{E}^S_{k},  e^\prime_{i} \geq e^\prime_{i_k}$. Since $S^*$ is optimal,  $B(S^*) = B(S^\prime)$ and from equation \Cref{trans_opt_sched}  $b_{i_k} = 0 $. 
\endproof

\begin{repeatcorollary}
 $\forall S \in \feasS\big([a_j]_{j \in \mathcal{A}}, W\big)$, the total potential bumps $B(S) = \sum\limits_{k\in A} b_{i_k}$.
\end{repeatcorollary}

This corollary directly follows from Lemma 1.

\begin{repeatlemma}

 $S(A) \in \feasS\big([a_j]_{j \in \mathcal{A}}, W\big)$. 
\end{repeatlemma}
\proof{\textbf{Proof of Lemma 2}} 
For $\ntpss$, any feasible schedule must satisfy the seniority constraint and the horizon constraint.  Any block schedule $S(A^\prime)$ already satisfies the seniority constraint. For the horizon constraint, the following equation must hold:

\begin{equation*}
    H = C^*_0  - W \geq   C^*_0  - \sum\limits_{j \in A^\prime } a_j = C(S(A))
\end{equation*}
Since $H \geq C(S(A))$, $S(A)$ is a feasible schedule.

\endproof

\begin{repeatlemma}

$\exists A^\prime \subseteq A$ such that  $S(A^\prime) \in \optS\big([a_j]_{j \in \mathcal{A}}, W\big)$. 
\end{repeatlemma}

\proof{\textbf{Proof of Lemma 3}} 
Let $S^* = [s^*_i, e^*_i]_{(i \in \mathcal{E})}$ be an optimal schedule with total potential bumps $B(S^*)$ and makespan $C(S^*) \leq H$ for $\ntpss$. From Lemma 2,  a feasible solution always exists. We are now going to show  $\forall S^* \in \optS\big([a_j]_{j \in \mathcal{A}}, W\big)$ there exists a  block schedule $S(A^\prime) = [s_i, e_i]_{(i \in \mathcal{E})}$ where $A^\prime \subseteq A$ such that $B(S^*) = B(S(A^\prime))$ and $C(S(A^\prime)) \leq C(S^*) \leq H$.\\
From Lemma 1, $\forall k \in A, b_{i_k} = a_k \text{ or } 0$. One can simply construct $A^\prime = \{k : b_{i_k} = a_k \}$.    
The block schedule $S(A^\prime)$ has total potential bumps $B(S(A^\prime)) = \sum\limits_{j \in A^\prime } a_j = B(S^*)$. 
Also $\forall k \in A$,  the following two cases exists:
\textbf{\\Case 1} $k \in A^\prime$,  $b_{i_k} = a_k$: For this condition to be true $\forall i \in \mathcal{E}^S_k, e_i < e_{i_k}$. Hence  $ s_{i_k}\leq s_i   \leq s_{i_{k+1}}$.\\
\textbf{Case 2}  $k \notin A^\prime$,  $b_{i_k} = 0$: For this condition to be true $\forall i \in \mathcal{E}^S_k, e_i \geq e_{i_k}$. This implies $s_i + r_i \geq s_{i_k} + r_{i_k}.$ Hence it follows that $ s_{i_k} + a_k \leq s_i \leq s_{i_{k+1}}$.\\
Combining the two cases to get:  
\begin{equation*}
    s_{i_N} \geq \sum\limits_{j=1}^{N-1}\mathbf{1}_{j \notin A^\prime} a_j 
\end{equation*}
Therefore, for  the last employee, $s_M \geq s_{i_N} + \mathbf{1}_{j \notin A^\prime} a_N \geq \sum\limits_{j \notin A^\prime}a_j $.
Hence $C(S(A^\prime)) = e_M + s_M + r_M \geq \sum\limits_{ j\notin A^\prime} a_j + \sum\limits_{ j\in A} a_j$ and $S(A^\prime)$ is also an optimal schedule. 
\endproof

\begin{repeatcorollary}
 $\forall S^* \in \optS\big([a_j]_{j \in \mathcal{A}}, W\big)$, the optimal objective value $B(S^*) \geq W$.
\end{repeatcorollary}

\proof{\textbf{Proof of Corollary 3}} 
From Lemma 3, $\exists A^\prime \subseteq A, S(A) \in \optS\big([a_j]_{j \in \mathcal{A}}, W\big)$ such that $B(S(A^\prime)) = \sum\limits_{j \in A^\prime} a_j$
and $C(S(A^\prime)) = C^*_0 - \sum\limits_{j \in A^\prime } a_j $.
\begin{equation*}
    C(S(A^\prime)) = C^*_0 - B^* \leq H = C^*_0 - W
\end{equation*}
Simplifying, $B^* \geq W$.
\endproof
\begin{repeattheorem}
For $\ntpss$, $ B(S^*) = W \Leftrightarrow \exists A^* \subseteq A $ such that $\sum_{j \in A^*} a_j = W$.
\end{repeattheorem}
\proof{\textbf{Proof of Theorem 2}} 
\textbf{\\Case 1} $\exists$  $S^* = [s^*_i, e^*_i]$ such that $B(S^*) = W$: From Lemma 3, one can build an optimal block schedule $S(A^*)$ with $B(S(A^*)) = \sum\limits_{j \in A^* } a_j$.  Hence, $W = \sum\limits_{j \in A^* } a_j$.\\
\textbf{Case 2}  Suppose  $ \exists A^* \subseteq A$ such that $W = \sum\limits_{j \in A^* } a_j$: Construct a feasible block schedule $S(A^*)$ with $B(S(A^*)) = \sum\limits_{j \in A^* } a_j = W$  and  makespan $C(S(A^*)) = C^*_0 - \sum\limits_{j \in A^* } a_j = H$. By Corollary 3, $S(A^*)$ is optimal.  
\endproof
\begin{repeattheorem}
 \textsc{NTP}$(I)$ is \emph{NP-complete}.
\end{repeattheorem}
\proof{\textbf{Proof of Theorem 3}} 
This conclusion follows from the fact that \textsc{Subset-Sum} is \emph{NP-complete}. It is straightforward to verify that \textsc{NTP}$(I)$ belongs to the class NP. Therefore, by providing a reduction from the \textsc{Subset-Sum} problem to \textsc{NTP}$(I)$, we establish the NP-completeness of \textsc{NTP}$(I)$.

\endproof

 \begin{repeattheorem}

Given any online algorithm, there exists an instance $I$ in which the algorithm suffers the maximum potential bumps or has vacant shifts, and the offline counterpart will have zero potential bumps and zero vacancies for that same instance.
\end{repeattheorem}
\proof{\textbf{Proof of Theorem 4}} 
Let $H$ be the horizon and $\mathcal{E}$ denote the employee set.  Let $s_i$ denote times when the online algorithm sends the notifications to employee $i$. The goal is to design a suitable instance (response delays $r_i$) such that the online algorithm behaves as worst as possible. In other words, we designing an adversary for the online algorithm. As the online algorithm starts sending notifications, two cases can happen

\textbf{\\Case 1} $\nexists$  $\hat{i}$ such that $\hat{i} = \{\min i : i \in \mathcal{E} , s_i >0\}$.\\
 Alternately this means all employees are sent notifications at time 0. In such a case set

  \begin{equation}
  {r_{i}} = H + 1 - i, \quad \forall i \in \mathcal{E}
 \end{equation}
The online algorithm will then suffer the maximum potential bumps as employees will respond in the reverse order of seniority ($e_i = H + 1 - i$). The offline algorithm, however, can easily prevent bumps by again simply ensuring that all employees respond back at the same time.  Let $s^*_i$ denote the offline solution.
  \begin{equation}
  {s_{i}^*} = H - r_i, \quad \forall i \in \mathcal{E}
 \end{equation}
 This will ensure all employees respond at the horizon.
\textbf{\\Case 2} $\exists$  $\hat{i}$ such that $\hat{i} = \{\min i : i \in \mathcal{E} , s_i >0\}$.\\
 $\hat{i}$ represents the first employee not notified at time = 0. In such a case, we can set 
  \begin{equation}
  r_{i} = H
 \end{equation}
 Since $e_i = s_i + H $, it would mean all employee junior employees to $\hat{i}$ and as well as $\hat{i}$ would respond after the horizon. Hence, the moment the online algorithm waits for even one unit we make sure that irrespective of the notification times for the next employees, we set their response times such that there are some vacant shifts.  
  It is easy to see that the offline solution will have no bumps due to $r_i = r_j, \forall i, j \in \mathcal{E}$. Additionally, we can also ensure an offline schedule has zero vacancies if 
 
\begin{equation}
\label{off_sched1}
s^*_{j} = 0
\end{equation}
This way, you can ensure that all employees respond at the horizon and that there are zero vacancies or bumps. We also respect seniority in the notifying order.

 \section{Full information (offline) formulation - ${\textbf{\textsc{MIP}}_{\textbf{NTP2}}}$}
This formulation explicitly considers shift vacancy cost instead of introducing a hard constraint to ensure complete shift scheduling. \Cref{tab:all_nots3} presents the parameters for the formulation.  The  formulation is similar to $\textbf{\textsc{MIP}}_{\textbf{DNTP-S}}$. However, since we are only concerned with one scenario at a time, we need only one single additional variable $z$ to track when an employee has responded. We have the same objective in \eqref{eqwo1} as in ${\textbf{\textsc{MIP}}_{\textbf{DNTP-S}}}$. Constraints \eqref{eqwo3} - \eqref{eqwo6} have the same meaning to Constraints \eqref{eqws3} - \eqref{eqws10} in ${\textbf{\textsc{MIP}}_{\textbf{DNTP-S}}}$. Given the high cost $G$ for a shift vacancy this formulation will always try to schedule all shifts first.
\begin{table}[ht!]
\parbox{1\linewidth}{
\centering
\resizebox{\linewidth}{!}{
  {\def\arraystretch{1.2}  
\begin{tabular}{llll}
\hline

\textbf{Parameters}&\\
\hline
$M$ &  number of employees\\
$L$ &  number of shifts\\
\hline
\textbf{Variables}&\\
\hline
$z_i$ & binary variable equal to 1 if employee $i \in \mathcal{E}$ responds within the horizon else 0\\
$\theta$ & number of shifts vacant \\ 
\hline
\end{tabular}}}

\caption{Parameters, and Variables of the {\dettwoprobname}}
\label{tab:all_nots3}
}
\end{table}

\begin{align} 
{\textbf{\textsc{MIP}}_{\textbf{NTP2}}} := \hspace{0.2cm }  \boldsymbol\min \quad & G \theta + \sum_{i\in \mathcal{E}}\sum_{j\in \mathcal{E}:i<j} y_{ij}  \label{eqwo1}
\end{align}
\textbf{Constraints}
\begin{align}
&s_i   \leq s_j & \forall i,j\in \mathcal{E}, i < j\label{eqwo2}\\
&s_i + r_i \geq  (H+1) (1 - z_i) & \forall i\in \mathcal{E} \label{eqwo3}\\
&s_i + r_i \leq  H + r_i (1- z_i) & \forall i\in \mathcal{E} \label{eqwo4}\\
&s_i - s_j + \delta_{ij} \leq \delta_{ij}y_{ij} +  (H + r_i)(1 - z_{i}) & \forall i,j\in \mathcal{E}, i < j,  r_j  \leq r_i \leq D \label{eqwo5}\\
&\sum\limits_{i\in \mathcal{E}}z_i   + \theta \geq L & \label{eqwo6}\\
&s_i + 1 \leq  s_{i + W} & \forall i, i + W \in \mathcal{E}\label{eqwo9}\\
&y_{ij} \in \{0,1\} & \forall i,j\in \mathcal{E}, i < j\label{eqwo7}\\
&s_{i} \in \mathbb{R}^+ & \forall i,\in \mathcal{E}\label{eqwo8}\\
&z_{i} \in \{0,1\} & \forall i\in \mathcal{E}\label{eqwo12}
\end{align}
The following additional constraints are added to the formulation. 

\begin{align}
&s_{i+1} - s_i \leq H (1 - z_i) & \forall i \in \mathcal{E} - \{M\}, r_{i+1} \geq r_i \label{eqwo13}\\
&y_{ij}  \leq z_i & \forall i,j\in \mathcal{E}, i < j\label{eqwo10}\\
&y_{ij}  \leq  z_j & \forall i,j\in \mathcal{E}, i < j\label{eqwo11}
\end{align}

Constraint \eqref{eqwo13} specifically refers to case when $r_{i+1} \geq r_i$.  Since the next employee, $i+1$, has a greater response time, there is no incentive to wait to notify $i+1$ unless that employee is to be pushed out of the horizon. Hence, Constraint \eqref{eqwo13} will ensure that this employee is either notified at the same time or responds after the end of the horizon. Constraints \eqref{eqwo10} and \eqref{eqwo11} ensures that employee $i$ can only bump employee $j$ if both of them respond within the time horizon. 

\section{Figures}

\begin{figure}[!htbp]
    \centering

    \begin{minipage}{.45\textwidth}
        \centering
        \caption{Pareto Front of all policies - Shift Vacancy \% vs Number of Bumps, $ D = 2$ }
\includegraphics[width=\textwidth]{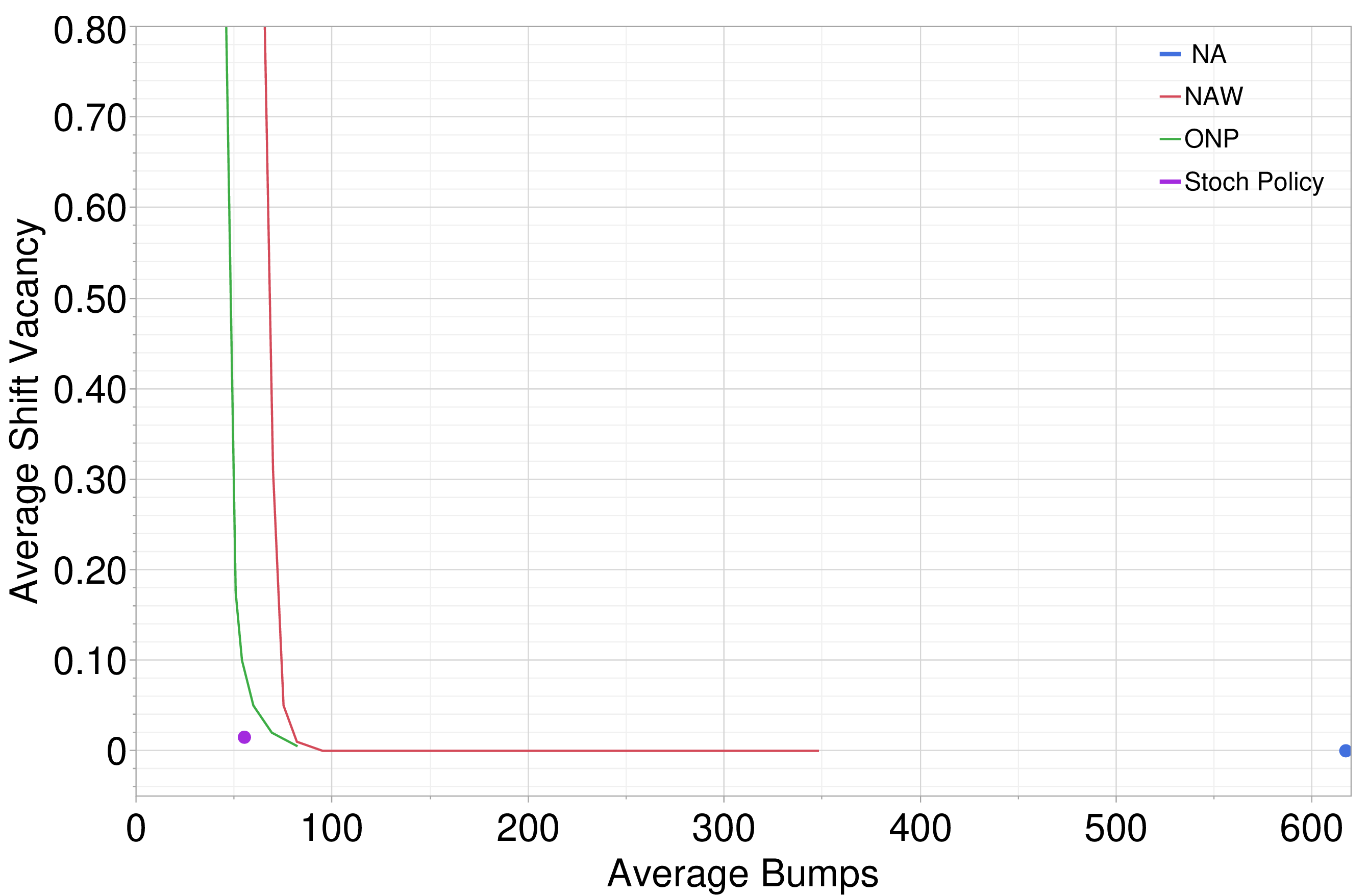}
    \label{resfigure1} 
     \end{minipage}%
    \begin{minipage}{.45\textwidth}
        \centering
    \caption{Pareto Front of all policies - Shift Vacancy \% vs Number of Bumps, $ D = 3$}
\includegraphics[width=\textwidth]{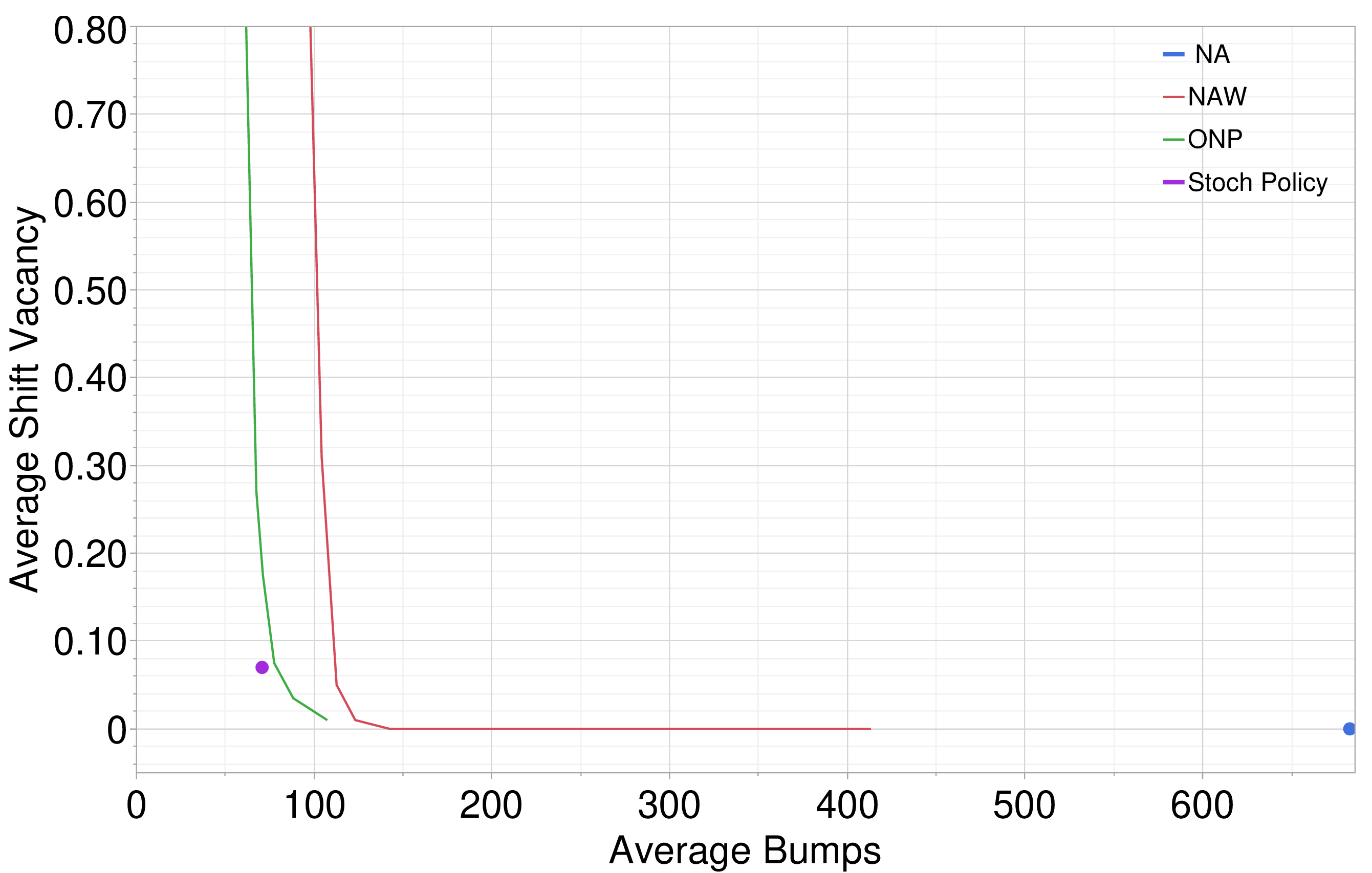}       \label{resfigure2}
    \end{minipage}
    
\end{figure}

\Cref{resfigure1}, \Cref{resfigure2}  compares all the policies for different values of $D$ on the 200 training instances that were used as scenarios to solve $\textbf{\textsc{MIP}}_{\textbf{DNTP-S}}$. The figures plot the Pareto frontier of the average number of bumps suffered by the employees against the average shifts vacant for ONP, NA, and NAW. The graphs contain one point for each unique parameter of the policy.  In the case of the dynamic policy based on offline solutions, this parameter would be the descriptive statistic used.   All policy points are connected by a line to form a Pareto frontier of that policy. Both plots show that the offline solution-based policy is consistently better than the heuristic policy. The stochastic solution obtained from $\textbf{\textsc{MIP}}_{\textbf{DNTP-S}}$ falls very close to the ONP frontier, suggesting we do very well with this policy. We also plot the performance of the Notify All (NA) policy. This policy will notify all employees in one go at the start of the horizon. We use this policy to demonstrate how bad the system can be if everyone is notified at once. As expected, we see a higher number of bumps when $D = 3 $ hours as compared to when $D = 2$ hours.